\documentclass[12pt]{amsart}
\usepackage[left=2.90cm, right=2.90cm, top=3.30cm, bottom=3.30cm]{geometry}

\usepackage{amsfonts,graphics,amsmath,amsthm,amsfonts,amscd,amssymb,amsmath,latexsym,multicol,
 mathrsfs}
\usepackage{epsfig,url}
\usepackage{flafter}
\usepackage{fancyhdr}
\usepackage[hyperindex,breaklinks]{hyperref}
\hypersetup{colorlinks=true, linkcolor=black, breaklinks=true}

\usepackage{graphicx}
\usepackage{xcolor}
\usepackage{float}
\usepackage{bm}
\usepackage{enumitem}
\usepackage{multirow}
\usepackage{scalerel,stackengine}
\usepackage{stmaryrd}
\usepackage{datetime}
\usepackage[title]{appendix}
\usepackage{url}
\usepackage{tabulary}
\usepackage{booktabs}
\usepackage{tikz}
\usetikzlibrary{positioning}
\usepackage{verbatim}

\usepackage[utf8]{inputenc}
\usepackage[english]{babel}
\usepackage{csquotes}

\usepackage{comment}

\usepackage{CJKutf8}

\usepackage{quiver}

\theoremstyle{definition}
\newtheorem*{Pf}{Proof}
\newtheorem*{Ack}{Acknowledgement}
\newtheorem*{Str}{Strategy}

\theoremstyle{plain}
\newtheorem{Thm}{Theorem}[section]
\newtheorem{Lem}[Thm]{Lemma}
\newtheorem{Rem}[Thm]{Remark}

\newtheorem{Def}[Thm]{Definition}
\newtheorem{Prop}[Thm]{Proposition}
\newtheorem{Cor}[Thm]{Corollary}

\newtheorem{Conj}[Thm]{Conjecture}

\title[\textbf{On finiteness of fiber space structures of klt CY pairs in dim 3}]{\textbf{On finiteness of fiber space structures of klt Calabi-Yau pairs in dimension 3}}
\author{Fulin Xu}
\date{(Dated: \today)}
\begin{document}

\begin{abstract}
We prove that for a fixed klt Calabi-Yau pair $(X,\Delta)$ of dimension $3$, the set of fiber space structures of $X$ is finite up to $Aut(X,\Delta)$. 
\end{abstract}

\maketitle

\tableofcontents

\section{Introduction}
Motivated by mirror symmetry, the Morrison-Kawamata cone conjecture, \cite{Mor93}, \cite{Kaw97}, predicts that the effective nef cone of a $K$-trivial variety admits a
rational polyhedral fundamental domain under the action of the automorphism group. An important geometric consequence is that the set of all possible contractions from a fixed $K$-trivial variety is finite up to automorphism. 

From the point of view of birational geometry, one should work on pairs with certain singularities rather than only work on smooth varieties. \cite{Tot09} generalized this conjecture for pairs: 

\begin{Conj}
Let $(X/S,\Delta)$ be a projective klt Calabi-Yau pair, where $X$ is $\mathbb{Q}$-factorial with $\Delta$ an $\mathbb{R}$-boundary. Then 

(a) The number of $Aut(X/S, \Delta)$-equivalence classes of faces of the cone $\mathcal{A}^e(X/S)$
corresponding to birational contractions or fiber space structures is finite. Moreover, the action of $Aut(X/S,\Delta)$ on the effective nef cone $\mathcal{A}^e(X/S)$ admits a rational polyhedral fundamental domain $\Pi$, in the sense that: 

(1) $\mathcal{A}^e(X/S) = \cup_{g\in Aut(X/S,\Delta)}g_*\Pi$, 

(2) $\text{Int }\Pi \cap g_*\text{Int }\Pi= \emptyset$ unless $g_* = 1$.

(b) The number of $PsAut(X/S, \Delta)$-equivalence classes of chambers $\mathcal{A}^e
(X_0 /S)$
in the cone $\mathcal{M}^e
(X/S)$ corresponding to marked small $\mathbb{Q}$-factorial modifications $X_0 \to S$ of $X \to S$ is finite. Moreover, the action of $PsAut(X/S,\Delta)$ on the effective movable cone $\mathcal{M}^e(X/S)$ admits a rational polyhedral fundamental domain $\Pi^\prime$, in the sense that: 

(1) $\mathcal{M}^e(X/S) = \cup_{g\in PsAut(X/S,\Delta)}g_*\Pi^\prime$, 

(2) $\text{Int }\Pi^\prime \cap g_*\text{Int }\Pi^\prime= \emptyset$ unless $g_* = 1$. 

\end{Conj}

\begin{Rem}
This will be called the cone conjecture for short. 
The first (resp. second) part of (a) is usually referred to as the weak nef cone conjecture (resp. nef cone conjecture), and the first (resp. second) part of (b) is usually referred to as the weak movable cone conjecture (resp. movable cone conjecture). Recall that a small $\mathbb{Q}$-factorial modification (SQM) of $X$ over $S$ means a pseudo-isomorphism over $S$ from $X$ to some other $\mathbb{Q}$-factorial variety with a projective
morphism to $S$.

In this paper, we mainly consider the absolute case, i.e. when $S$ is a single point. 

Assuming standard conjectures in MMP, existence of a rational polyhedral fundamental domain for $\mathcal{M}^e(X)$ implies all the other statements by \cite[Section 2.5]{Xu24}. 
It is worth noting that \cite{GLSW24} gives more equivalent formulations up to standard conjectures in MMP. 

Also, it is known that it makes no difference to assume $\Delta$ to be a $\mathbb{Q}$-divisor. Indeed, we may approximate $\Delta$ by $\mathbb{Q}$-divisors $\Delta^\prime$ such that $\text{Supp} \Delta = \text{Supp} \Delta^\prime$, $K_X+\Delta^\prime \equiv 0$, and the differences of the coefficients are arbitrarily small. 
\end{Rem}

In dimension $2$, the cone conjecture is solved in full generality by \cite{Tot09}.

In higher dimensions, it is hard to deal with the general case. However, the Beauville-Bogomolov decomposition, \cite{Bea83}, shows that, up to an \'etale
cover, every smooth projective $K$-trivial variety can be decomposed into a product of abelian, hyperk\"ahler, and
strict Calabi–Yau manifolds. We note that a singular analog is established in \cite{HP17}. 

For abelian varieties, the cone conjecture is solved by \cite{PS12}. For hyperk\"ahler varieties, the cone conjecture is solved to a large extent, see \cite{Mar11} for a good survey. However, the cone conjecture is widely open for strict Calabi-Yau varieties. Also, it is worth noting that the cone conjecture is sensitive with respect to taking \'etale cover, so solutions in these three cases do not imply the general case immediately. \cite{PS23} proved the cone conjecture for certain Enriques manifolds. \cite{MQ24} proved the cone conjecture for generalized hyperelliptic varieties. 

We briefly summarize some known general facts about the cone conjecture of strict Calabi-Yau manifolds: 

In dimension $3$, we have a series of partial results. \cite{OP98} proved that every strict Calabi-Yau threefold with positive second Chern class admits only finitely many fiber space structures. \cite{Sze99} proved finiteness of ample divisors of bounded volume up to automorphism. 
\cite{OS01} proved finiteness of $c_2$-contractions up to automorphism, and in particular finiteness of abelian fibrations up to automorphism. 
\cite{Ueh04} proved finiteness of divisorial contractions whose exceptional divisor has trivial intersection with $c_2$. 
Recently, \cite{Lut24} proved that the cone conjecture is preserved under deformations. 

In general dimension, if the Picard number is small, one can still have some control. See \cite{LP13}, \cite{Ogu14} for the Picard number $2$ case. 

We note that there are also many interesting results for the cone conjecture in the relative setting, as well as various results for many explicit examples. We also recommend 
\cite{LOP18} for a good survey.

The following is our main result: 

\begin{Thm}\label{Thm main}
Let $(X,\Delta)$ be a projective $\mathbb{Q}$-factorial klt Calabi-Yau pair in dimension $3$. Then

(1) If $\kappa(-K_X)>0$, then the cone conjecture holds for $(X,\Delta)$. 

(2) If $\kappa(-K_X)=0$, we may run $K_X+(1+\varepsilon)\Delta$-MMP for a small positive rational number $\varepsilon$ and terminate with $X^{min}$. Let $\tilde{q}(X^{min})$ be the augmented irregularity of $X^{min}$. Then 

(a) If $\tilde{q}(X^{min})>0$, then the cone conjecture holds for $(X,\Delta)$. 

(b) If $\tilde{q}(X^{min})=0$, then there are only finitely many fiber space structures of $X$ up to $Aut(X,\Delta)$. 
\end{Thm}

\begin{Rem}
Part (1) is due to \cite{Xu24}, and we include this for completeness of the statement. 

The augmented irregularity $\tilde{q}(X)$ is defined to be the maximum of $h^1(\tilde{X},\mathcal{O}_{\tilde{X}})$, where $\tilde{X} \to X$ is a finite quasi-\'etale cover. 
\end{Rem}

When $\Delta = 0$, the statements become much simpler: 

\begin{Cor}
Let $X$ be a projective $\mathbb{Q}$-factorial klt variety in dimension $3$ with $K_X\equiv 0$. Let $\tilde{q}(X)$ be the augmented irregularity of $X$. Then 

(a) If $\tilde{q}(X)>0$, then the cone conjecture holds for $X$. 

(b) If $\tilde{q}(X)=0$, then there are only finitely many fiber space structures of $X$ up to $Aut(X)$. 
\end{Cor}

\begin{Rem}
Our result gives new positive evidence for the cone conjecture in dimension 3. We essentially reduce the cone conjecture in dimension 3 to the strict Calabi-Yau case, and prove finiteness of fiber space structures. 

However, it seems that we are still far from proving the cone conjecture (even finiteness of birational contractions) in dimension 3. 
\end{Rem}

\begin{Rem}
For case (2)(a) of the main theorem, a typical example is the following: 

Consider an algebraic K3 surface $X$ with an involution $\tau$ acting on $H^0(X,\omega_X)$ by $-1$, and an elliptic curve $E$ with an involution $\iota$ acting by multiplication of $-1$. Then $Y = E\times X/\iota\times \tau$ is a canonical Calabi-Yau threefold with positive augmented irregularity. 

If one takes a crepant resolution $\tilde{Y}$ of $Y$, then the resulting threefold is studied by Borcea and Voisin \cite{Bor97}. However, it seems still hard to prove the cone conjecture for the crepant resolution $\tilde{Y}$.

\end{Rem}

\begin{Cor}\label{Cor main}
Let $(X,\Delta)$ be a projective klt Calabi-Yau pair in dimension $3$. Then $X$ admits only finitely many fiber space structures up to $Aut(X,\Delta)$. 
\end{Cor}

\begin{Str}
We use quite different arguments for part (a) and part (b) of Theorem \ref{Thm main}(2). We will call them non-strict case and strict case respectively. 

For the strict case, our strategy is to deduce finiteness from boundedness and discreteness. 

Let us start with a motivating example. For elliptic fibrations of strict Calabi-Yau threefolds, we already have enough tools. Boundedness of elliptic Calabi-Yau threefolds is established in \cite{FHS21}, and discreteness is essentially done in \cite{Sze99}. So the proof would just be a combination of these inputs. 

For general curve fibrations, we need to overcome more difficulties. For boundedness, we need a generalization of \cite{FHS21} for a certain type of pairs. For discreteness, we need to generalize some arguments in \cite{Sze99} on period maps to the singular case. Fortunately, these generalizations work for a suitably chosen class of Calabi-Yau pairs. 

For surface fibrations, an interesting observation is that if one has at least two different surface fibrations, then they are dominated by a curve fibration. In order to prove finiteness, we may assume all surface fibrations factor through some curve fibration, say $X \stackrel{\pi}{\to} B \stackrel{p}{\to} C$. As the cone conjecture is known for klt surface pairs by \cite{Tot09}, the $\pi$ and $p$ are finite up to $Aut(X)$ and $Aut(B)$, respectively. Now we need a lifting property from $Aut(B)$ to $Aut(X)$. It is interesting that this lifting property can be deduced from a generalization of the boundedness result of curve fibered Calabi-Yau pairs. 

For the non-strict case, our strategy is to decompose the pair $(X,\Delta)$ via the Albanese map. 

For klt Calabi-Yau pairs, the Albanese map is an isotrivial fibration by \cite{Kaw85}, \cite{Amb16}. 
After a finite cover, $(X,\Delta)$ is either an abelian variety or a product of a curve and a surface. 
In the former case, the cone conjecture is essentially done by \cite{MQ24}. 
In the latter case, we may reduce to the lower dimensional cases. 
\end{Str}

\begin{Ack}
The author would like to thank Professor Caucher Birkar for his constant support and guidance. He thanks Professor Keiji Oguiso for kind suggestions. He thanks Professor Tsung-Ju Lee, Professor Mao Sheng, and Professor Dingxin Zhang for expertise in Hodge theory. He thanks Xintong Jiang for reading a draft of this paper and giving useful comments. He thanks Yixuan Fu, Xiaowei Jiang, Xintong Jiang, Junpeng Jiao, Long Wang, and Zekai Yu for useful discussions. 
\end{Ack}

\section{Preliminaries}

Throughout this paper, we work over $\mathbb{C}$. 

Let $f:X \to S$ be a projective surjective
morphism of normal varieties, with geometrically connected fibers. Such an $f$ is called a contraction. If $f$ is birational, then $f$ is called a birational contraction. If $f$ is not birational, 
$f$ is called a fibration or a fiber space structure. 

\subsection{Types of divisors}

We recall some standard definitions for different kinds of Cartier divisors. 

Let $D$ be a Cartier divisor on $X$. It is said to be $f$-nef if $(D . C) \geq 0$ holds for any curve $C$ on $X$ which is
contracted by $f$. It is said to be $f$-movable if $\operatorname{codim} \text{Supp } \text{Coker}(f^*f_*\mathcal{O}_X(D) \to \mathcal{O}_X(D)) \geq 2$. It is said to be $f$-effective if $f_*\mathcal{O}_X(D) \neq 0$. And it is said to be $f$-big if its restriction to the generic fiber has maximal Kodaira dimension. 

\subsection{Various cones}

The N\'eron-Severi group
$$NS
(X/S) = \{\text{Cartier divisor on }X\}/(\text{numerical equivalence over }S)$$
is known to be finitely generated. Let $N^1
(X/S) = NS
(X/S)\otimes_{\mathbb{Z}} \mathbb{R}$. We set $\rho(X/S) = \dim N^1
(X/S)$ to be the relative Picard number. The numerical class of an $\mathbb{R}$-Cartier divisor $D$ in $N^1(X/S)$ is denoted by $[D]$.

The $f$-nef cone is denoted by $\bar{\mathcal{A}}(X/S)$, the closed $f$-movable cone is denoted by $\bar{\mathcal{M}}(X/S)$, and the $f$-pseudo-effective cone is denoted by $\bar{\mathcal{B}}(X/S)$. 
They are defined to be the closed convex cones in $N^1
(X/S)$ generated
by the numerical classes of $f$-nef divisors, $f$-movable divisors and $f$-effective
divisors respectively. 

We have the following obvious inclusion relations:
$$\bar{\mathcal{A}}(X/S) \subset \bar{\mathcal{M}}(X/S) \subset \bar{\mathcal{B}}(X/S) \subset N^1
(X/S)$$

The interior $\mathcal{A}(X/S) \subset \bar{\mathcal{A}}(X/S)$  (resp. $\mathcal{B}(X/S) \subset \bar{\mathcal{B}}(X/S)$) is the open convex subcone generated by the numerical classes of $f$-ample divisors ($f$-big divisors) and it is called
an $f$-ample cone (resp. $f$-big cone). Let $\mathcal{M}(X/S)\subset \bar{\mathcal{M}}(X/S)$ be the subcone generated by $f$-movable divisors. 

We denote the $f$-effective cone by $\mathcal{B}^e 
(X/S)$, i.e. the convex cone
generated by $f$-effective Cartier divisors. We call $\mathcal{A}^e
(X/S) = \bar{\mathcal{A}}(X/S) \cap \mathcal{B}^e
(X/S)$
and $\mathcal{M}^e
(X/S) = \bar{\mathcal{M}}(X/S) \cap \mathcal{B}^e
(X/S)$ the $f$-effective $f$-nef cone and $f$-effective
$f$-movable cone, respectively.

We collect some lemmas, which will be used in the proof of our main results: 

\begin{Lem}\label{Lem pullback of divisors}
Let $f: X\to Y$ be a finite surjective morphism between normal projective varieties. Let $D$ be a divisor on $Y$. Then 

(1) If $D$ is linearly equivalent to some effective divisor, so is $f^*D$. Conversely, if $f^*D$ is linearly equivalent to some effective divisor, so is $nD$ for some $n\in \mathbb{N}_{>0}$. 

(2) $D$ is ample iff $f^*D$ is; 

(3) $D$ is nef iff $f^*D$ is; 

(4) If $D$ is movable, so is $f^*D$. Conversely, if $f^*D$ is movable, so is $nD$ for some $n\in \mathbb{N}_{>0}$. 

\end{Lem}

\begin{Pf}
We first prove (1). It is clear that $D$ is linearly equivalent to some effective divisor implies $f^*D$ is linearly equivalent to some effective divisor. Conversely, if $f^*D$ is effective, then by \cite[Tags 0BCY and 0BD3]{stacks-project}, we may take the norm of $\mathcal{O}_X(f^*D)$, say $\text{Norm}_f \mathcal{O}_X(f^*D)$, which is a line bundle on $Y$. Let $n$ be the degree of $f$. Moreover, by \cite[Tag 0BCY]{stacks-project}, $\text{Norm}_f \mathcal{O}_X(f^*D) \simeq \mathcal{O}_Y(nD)$, and by \cite[Tag 0BCZ]{stacks-project} with $\mathcal{L} = \mathcal{O}_X$ and $\mathcal{L}^\prime = \mathcal{O}_X(f^*D)$, we conclude that $nD$ is effective. 

Assertion (2) follows from \cite[Tag 0BD0]{stacks-project}. 

To prove (3), take an ample divisor $H$ on $Y$, then $f^*H$ is ample by (2). Then the result follows from the fact that $D$ is nef iff $nD+H$ is nef for any $n \in \mathbb{N}_{>0}$. 

It remains to show (4). Suppose $D$ is movable, then for any prime divisor $S$ on $Y$, there exists some $0 \leq E \sim D$ such that $S \not\subset E$. Let $T$ be any prime divisor on $X$, then take $S = f(T)$, and we have $0 \leq f^*E \sim f^*D$, $T \not\subset f^*E$. So $f^*D$ is movable. 

Conversely, suppose $f^*D$ is movable. Possibly replacing $X$ by the Galois closure, we may assume $f$ is Galois with Galois group $G$, and $f^*D$ remains movable. Let $n$ be the degree of the Galois morphism $f$. For any prime divisor $S$ on $Y$, we need to exhibit $0 \leq E \sim nD$ such that $S \not\subset E$. Indeed, for any prime divisor $T \subset f^*S$, by assumption, there exists a rational function $\varphi_T$ on $X$ such that $0 \leq E_T = (\varphi_T)+f^*D \sim f^*D$ satisfies $T \not\subset F$. Let $\varphi$ be a general linear combination of $\{\varphi_T\}_{T\subset f^*S}$, then $0 \leq E = (\varphi)+f^*D \sim f^*D$ satisfies $T \not\subset E$ for any prime divisor $T\subset f^*S$. 

Let $\psi$ be the product of Galois conjugations of $\varphi$, we claim that $0 \leq F = (\psi)+nD \sim nD$ satisfies $S \not\subset E$. Indeed, we have $f^*F = \sum_{g\in G} (g^*\varphi)+nf^*D = \sum_{g\in G} g^*E$, which is effective and does not contain any component of $f^*S$. So our claim holds. 

\hfill
\qedsymbol
\end{Pf}

A finite Galois morphism $f:X\to Y$ between normal varieties means a finite surjective morphism such that $\mathbb{C}(X)/\mathbb{C}(Y)$ is a Galois extension of fields. 
The variety $X$ is the normalization of $Y$ in $\mathbb{C}(X)$. 
Every element of $G=\operatorname{Gal}(\mathbb{C}(X)/\mathbb{C}(Y))$ preserves the integral closure of each affine coordinate ring of $Y$, and hence extends uniquely to an automorphism of $X$ over $Y$.

\begin{Lem} \label{Lem cone of finite map}
Let $f: X\to Y$ be a finite Galois morphism between normal projective varieties with Galois group $G$. Then the pullback map $f^*$ induces isomorphisms

(1) $N^1(Y) \simeq N^1(X)^G$

(2) $\mathcal{A}(Y) \simeq \mathcal{A}(X)^G$

(3)  $\mathcal{B}^e(Y) \simeq \mathcal{B}^e(X)^G$

(4) $\bar{\mathcal{A}}(Y) \simeq \bar{\mathcal{A}}(X)^G$

(5) $\bar{\mathcal{A}}(Y)_+ \simeq (\bar{\mathcal{A}}(X)_+)^G \simeq (\bar{\mathcal{A}}(X)^G)_+$

(6) $\mathcal{A}^e(Y) \simeq \mathcal{A}^e(X)^G$

(7) $\mathcal{M}(Y) \simeq \mathcal{M}(X)^G$

(8) $\mathcal{M}^e(Y) \simeq \mathcal{M}^e(X)^G$

\end{Lem}

\begin{Pf}
Let $d = |G| = \text{deg} f$. 

(1) We have injection $f^*: N^1(Y) \to N^1(X)$, whose image clearly lies in $N^1(X)^G$. On the other hand, let $\mathcal{L}$ be a line bundle on $X$ whose numerical equivalence class is fixed by $G$. Then $\mathcal{M}: = \otimes_{g\in G} g^*\mathcal{L}$ is numerically equivalent to $\mathcal{L}^{\otimes d}$ and $\mathcal{M}$ is a $G$-invariant line bundle. By \cite[Tags 0BCY and 0BD3]{stacks-project}, we may take the norm of $\mathcal{M}$, say $\text{Norm}_f \mathcal{M}$, which is a line bundle on $Y$. Moreover, by comparing the transition functions, we have $f^*\text{Norm}_f \mathcal{M} \simeq \mathcal{M}^{\otimes d}$. So surjectivity holds. 

(2) This follows from (1) and Lemma \ref{Lem pullback of divisors}(2), since ampleness is invariant under numerical equivalence. 

(3) It follows from (1) and Lemma \ref{Lem pullback of divisors}(1) that the natural map $\mathcal{B}^e(Y) \to \mathcal{B}^e(X)^G$ is well-defined and injective. To show surjectivity, we need to show that if $f^*\mathcal{L}$ is numerically equivalent to an effective line bundle, then it is numerically equivalent to a pullback of an effective line bundle. 

Suppose $f^*\mathcal{L}$ is numerically equivalent to an effective line bundle $\mathcal{N}$. Then $\mathcal{M}: = \otimes_{g\in G} g^*\mathcal{N}$ is numerically equivalent to $\mathcal{N}^{\otimes d}$ and $\mathcal{M}$ is a $G$-invariant effective line bundle. By \cite[Tags 0BCY and 0BD3]{stacks-project}, we may take the norm of $\mathcal{M}$, say $\text{Norm}_f \mathcal{M}$, which is a line bundle on $X$. Moreover, by comparing the transition functions, we have $f^*\text{Norm}_f \mathcal{M} \simeq \mathcal{M}^{\otimes d}$. Furthermore, applying \cite[Tag 0BCZ]{stacks-project} with $\mathcal{L} = \mathcal{O}_X$ and $\mathcal{L}^\prime = \mathcal{M}$, we conclude that $\text{Norm}_f \mathcal{M}$ is effective. So surjectivity holds. 

(4) This follows from (2) and the fact that $\bar{\mathcal{A}}(X)^G$ is closed and the relative interior is $\mathcal{A}(X)^G$. 

(5) This follows from (4) and the fact that the $G$ action on $N^1(X)$ respects $\mathbb{Q}$-structure. 

(6) This follows from taking intersection of (3) and (4). 

(7) The proof is essentially the same as (3). It follows form (1) and Lemma \ref{Lem pullback of divisors}(4) that the natural map $\mathcal{M}(Y) \simeq \mathcal{M}(X)^G$ is well-defined and injective. To show surjectivity, we need to show that if $f^*\mathcal{L}$ is numerically equivalent to a movable line bundle, so is $\mathcal{L}$. 

Suppose $f^*\mathcal{L}$ is numerically equivalent to a movable line bundle $\mathcal{N}$. Then $\mathcal{M}: = \otimes_{g\in G} g^*\mathcal{N}$ is numerically equivalent to $\mathcal{N}^{\otimes d}$ and $\mathcal{M}$ is a $G$-invariant movable line bundle. By \cite[Tags 0BCY and 0BD3]{stacks-project}, we may take the norm of $\mathcal{M}$, say $\text{Norm}_f \mathcal{M}$, which is a line bundle on $X$. Moreover, by comparing the transition functions, we have $f^*\text{Norm}_f \mathcal{M} \simeq \mathcal{M}^{\otimes d}$. Furthermore, by \cite[Tag 0BCZ]{stacks-project} with $\mathcal{L} = \mathcal{O}_X$ and $\mathcal{L}^\prime = \mathcal{M}$, we conclude that $\text{Norm}_f \mathcal{M}$ is movable, since the base locus must come from base locus of $\mathcal{M}$. So surjectivity holds. 

(8) This follows from taking closure of (7) and then taking intersection with (3). 

\hfill
\qedsymbol
\end{Pf}

\begin{Prop}\label{Prop 52 ex}
Let $X$ be a normal projective variety such that $H^1(X,\mathcal{O}_X) = 0$. Let $T$ be another variety. 

(a) If $\mathcal{L}$ is an invertible sheaf on $X \times T$, then the invertible sheaves $\mathcal{L}_t$, on 
$X = X \times \{t\}$ are isomorphic, for all closed points $t\in T$. 

(b) $Pic(X \times T) = Pic (X) \times Pic (T)$. 
\end{Prop}

\begin{Pf}
This follows from \cite[Chapter III, Exercise 12.6]{Har77}. 

\hfill
\qedsymbol
\end{Pf}

\begin{Cor}\label{Cor cone of product}
Let $X$, $Y$ be normal projective varieties such that $H^1(X,\mathcal{O}_X) = 0$. Then 

(1) $N^1(X\times Y) \simeq N^1(X)\times N^1(Y)$

(2) $\mathcal{B}^e(X\times Y) \simeq \mathcal{B}^e(X)\times \mathcal{B}^e(Y)$

(3) $\bar{\mathcal{A}}(X\times Y) \simeq \bar{\mathcal{A}}(X)\times \bar{\mathcal{A}}(Y)$. 

(4) $\mathcal{A}^e(X\times Y) \simeq \mathcal{A}^e(X)\times \mathcal{A}^e(Y)$. 

(5) $\mathcal{M}^e(X\times Y) \simeq \mathcal{M}^e(X)\times \mathcal{M}^e(Y)$. 
\end{Cor}

\begin{Pf}
Let $p_X$, $p_Y$ be projections of $X\times Y$ to $X$ and $Y$. Let $x$ (resp. $y$) be a closed point on $X$ (resp. $Y$). 

(1) We have a natural map $f: N^1(X)\times N^1(Y) \to N^1(X\times Y)$ given by pullback via two natural projections. Consider $N^1(X\times Y) \to N^1(X)\times N^1(Y)$, $\mathcal{L} \mapsto (\mathcal{L}|_{X\times\{y\}},\mathcal{L}|_{\{x\}\times Y})$, which is a left inverse of $f$. So $f$ is injective. By proposition \ref{Prop 52 ex}, $f$ is surjective, so the result holds. 

(2) This follows directly from the K\"unneth formula, see \cite[tag 0BED]{stacks-project}. Indeed, $H^0(X\times Y, p_X^*\mathcal{L}_X \otimes p_Y^*\mathcal{L}_Y) =H^0(X, \mathcal{L}_X) \otimes H^0(Y, \mathcal{L}_Y)$. 

(3) Clearly we have a natural injective map $f: \bar{\mathcal{A}}(X)\times \bar{\mathcal{A}}(Y) \to \bar{\mathcal{A}}(X\times Y)$, since pullbacks of nef divisors are nef. Conversely, if $p_X^*\mathcal{L}_X \otimes p_Y^*\mathcal{L}_Y$ is nef, by taking restriction to $X\times\{y\}$ and $\{x\}\times Y$, we can see both $\mathcal{L}_X$ and $\mathcal{L}_Y$ are nef. So $f$ is also surjective. 

(4) This is a direct corollary of (2) and (3). 

(5) Again by the K\"unneth formula, see \cite[tag 0BED]{stacks-project}, we have $H^0(X\times Y, p_X^*\mathcal{L}_X \otimes p_Y^*\mathcal{L}_Y) =H^0(X, \mathcal{L}_X) \otimes H^0(Y, \mathcal{L}_Y)$. In particular, we have ${\rm Bs}|p_X^*\mathcal{L}_X \otimes p_Y^*\mathcal{L}_Y| = {\rm Bs}|\mathcal{L}_X| \times Y\cup X \times {\rm Bs}|\mathcal{L}_Y|$. So $p_X^*\mathcal{L}_X \otimes p_Y^*\mathcal{L}_Y$ is movable if and only if both $\mathcal{L}_X$ and $\mathcal{L}_Y$ are movable. 

\hfill
\qedsymbol
\end{Pf}

\subsection{Pairs} 

We also recall some notations for pairs. 

For a normal variety $X$, let $K_X$ be the canonical divisor, $\omega_X$ be the canonical sheaf. When $X$ is projective, the canonical sheaf coincides with the dualizing sheaf by \cite[Proposition 5.75]{KM98}. Moreover, if $X$ is klt, $\omega_X$ coincides with the dualizing complex. 

For a reflexive sheaf $\mathcal{F}$, we denote the reflexive hull of the $n$-th tensor power of $\mathcal{F}$ by $\mathcal{F}^{[n]}$. We will only use the case when $\mathcal{F} = \omega_X$. 

A pair is a normal variety $X$ together with an $\mathbb{R}$-divisor $\Delta = \sum_i a_i\Delta_i$ such that $K_X+\Delta$ is $\mathbb{R}$-Cartier and $a_i\in [0,1]$. We adopt the usual definitions of singularities of pairs (lc, klt, canonical, terminal). 

Let $(X,\Delta)$ be a pair and let $f:X\to S$ be a morphism. We write $(X/S,\Delta)$ for these data. For an arbitrary morphism $T\to S$, the scheme-theoretic fiber product $X_T:=X\times_S T$ is defined, but need not underlie a pair. We use the notation $(X_T,\Delta_T)$ only when $X_T$ is normal, every $\Delta_i\times_S T$ is either empty or of pure codimension one in $X_T$, and the associated cycle \[ \Delta_T:=\sum_i a_i\Delta_{T,i},\qquad \Delta_{T,i}:=[\Delta_i\times_S T]_{\mathrm{codim}\,1}, \] is a boundary with $K_{X_T}+\Delta_T$ $\mathbb{R}$-Cartier. 
In applications, 
when $T \to S$ is \'etale, these conditions are automatically satisfied; 
when $f$ is a family and $T$ is a closed point on $S$, we need to restrict or stratify the base so that the normality and divisor conditions hold.

We say a pair $(X/S,\Delta)$ is Calabi-Yau if $K_X+\Delta \equiv_S 0$. By \cite[Corollary 1.6]{HX11}, we automatically have $K_X+\Delta \sim_{\mathbb{R},S} 0$. 

Let $(X,\Delta_X)$, $(Y,\Delta_Y)$ be two pairs. An isomorphism from $(X,\Delta_X)$ to $(Y,\Delta_Y)$ is an isomorphism $\phi: X\to Y$ such that $\phi_* \Delta_X= \Delta_Y$. A pseudo-isomorphism from $(X,\Delta_X)$ to $(Y,\Delta_Y)$ is a small birational map  $\phi: X\dashrightarrow Y$ such that $\phi_* \Delta_X= \Delta_Y$. In a similar way we define automorphisms and pseudo-automorphisms.

We denote the abstract group of automorphisms of the pair by
$Aut(X,\Delta)$ and the abstract group of pseudo-automorphisms by
$PsAut(X,\Delta)$. When $X$ is projective, we write
$\mathcal{A}ut(X,\Delta)$ for the automorphism group scheme over
$\mathbb{C}$ preserving the boundary, including its coefficients.
This group scheme is locally of finite type, and our convention is
\[
  Aut(X,\Delta)=\mathcal{A}ut(X,\Delta)(\mathbb{C}).
\]
We write $\mathcal{A}ut^0(X,\Delta)$ for its identity component and set
\[
  Aut^0(X,\Delta)=\mathcal{A}ut^0(X,\Delta)(\mathbb{C}).
\]
Thus $Aut^0(X,\Delta)$ denotes an abstract subgroup of
$Aut(X,\Delta)$. We use the same conventions for $Aut(X)$ and
$Aut^0(X)$ when the boundary is omitted, and for automorphisms over
a fixed base. Statements about morphisms, connected components, or
finite type refer to the group schemes; group actions and indices
refer to their groups of $\mathbb{C}$-points.

\subsection{Index 1 covers}

First we prove an elementary lemma. 

\begin{Lem}\label{Lem elementary}
Let $q = k/n$ be a rational number, $0<q<1$, $n, k \in \mathbb{N}_{>0}$, $(k,n) = 1$. Then the following conditions are equivalent: 

(1) $k = n-1$

(2) For any $i \in \mathbb{Z}$, $\lfloor (1-i) q\rfloor+ \lfloor i q\rfloor = 0$
\end{Lem}

\begin{Pf}
(1) $\Rightarrow$ (2): $-1< - \frac{n-1}{n} = (1-i) q - \frac{n-1}{n}+ i q - \frac{n-1}{n} \leq \lfloor (1-i) q\rfloor+ \lfloor i q\rfloor \leq (1-i) q + i q <1$. So $\lfloor (1-i) q\rfloor+ \lfloor i q\rfloor = 0$ since it is integral. 

(2) $\Rightarrow$ (1): Choose $i$ such that $ik \equiv n-1$ mod $n$. Such $i$ exists since $k$ and $n$ are coprime. Then $\lfloor (1-i) q\rfloor+ \lfloor i q\rfloor = -1$ unless $k = n-1$. 

\hfill
\qedsymbol
\end{Pf}

We use the following global index-one-cover construction. The standard
coefficients $1-1/n$ ensure that the ramification divisor cancels the
pullback of the boundary, so that the cover has no residual boundary.
Equivalently, they give the rounding identity in Lemma \ref{Lem elementary}.
The local case (1)(2)(3) are given by \cite{Ko13}. Although the proof of the global case is parallel, it is hard to find a direct reference.  The lifting assertion in (4) is extracted from
\cite[Lemma 39]{Xu24}, with simplifications suggested by Xintong Jiang:

\begin{Lem}\label{Lemma delicate lifting property}
Let $(X,\Delta)$ be a klt pair such that $K_X + \Delta$ is $\mathbb{Q}$-linearly trivial and the
coefficients of $\Delta$ are of the form $1-1/n$ for some integer $n$. 

Fix a section 
$\sigma \in H^0(X,\omega_X^{[m]}(m\Delta))$
, where $m$ is the global index, i.e. the minimal positive integer $n$ such that $n(K_X+\Delta)$ is Cartier and linearly equivalent to $0$. Consider the index $1$ cover of $(X,\Delta)$ corresponding to $\sigma$, say $p: Y \to X$, as defined in \cite[Definition 2.49]{Ko13}, i.e. $Y$ is the normalization of $ Spec_{X}\oplus_{i=0}^{m-1} \omega_X^{[i]}(\lfloor i\Delta \rfloor)$, where the multiplication structure is given by $\omega_X^{[m]}(m\Delta) \to \mathcal{O}_X$, $s\mapsto s\sigma^{-1}$. 

Then we have: 

(1) $Y$ is integral, 

(2) $Y = Spec_{X}\oplus_{i=0}^{m-1} \omega_X^{[i]}(\lfloor i\Delta \rfloor)$, 

(3) $K_Y \sim 0$, and $Y$ has canonical singularities. 

(4) Suppose that $X$ is proper, then any automorphism of $(X,\Delta)$ lifts to an automorphism of $Y$. 
\end{Lem}

\begin{Pf}

(1) 
Let $s$ be any rational section of $\omega_X$, then $s^m$ is a rational section of $\omega_X^{[m]}(m\Delta)$. So we have a function $\sigma/s^m \in \mathbb{C}(X)$. To show $Y$ is integral, it suffices to restrict to the generic fiber, i.e.
it suffices to show the polynomial $q(x) = x^m - \sigma/s^m$ is irreducible in $\mathbb{C}(X)$. Let $x_0$ be a root of $q(x)$, then $q(x) = \prod_{i=0}^{m-1} (x-x_0\xi^i)$, where $\xi$ is an $m$-th root of unity. Suppose $q(x)$ is reducible in $\mathbb{C}(X)$, then there exists some $0<d<m$ such that $x_0^d \in \mathbb{C}(X)$. Let $d_0$ be the greatest common divisor of $m$ and $d$, then there exist some integers $a$, $b$ such that $ad+bm=d_0$. So $x_0^{d_0} = (x_0^d)^a(x_0^n)^b \in \mathbb{C}(X)$. We may assume $d|m$ by replacing $d$ by $d_0$. Let $r = m/d$. Consider $t = x_0^d s^d$, which is a rational section of $\omega_X^{[d]}(\lfloor d\Delta \rfloor)$. Moreover, since $t^r = \sigma \in H^0(X,\omega_X^{[m]}(m\Delta))$, $d\Delta$ must be integral and 
$t$ gives a trivialization of $\omega_X^{[d]}(d\Delta)$. This contradicts the minimality of $m$.

(2) We follow the arguments in \cite[2.44]{Ko13}. 

It suffices to check that the right hand side is normal, which can be checked locally on $X$. So we may assume $X$ is affine. Then for any coherent $\mathcal{O}_X$-module $\mathcal{F}$ on $X$, we may use the same notation to denote the corresponding $H^0(X,\mathcal{O}_X)$-module $H^0(X,\mathcal{F})$.

Let $g$ be the generator of $Gal(\mathbb{C}(Y)/\mathbb{C}(X))$ that acts on $\mathbb{C}(X)\omega_X^{[i]} \subset \mathbb{C}(X)\oplus_{i=0}^{m-1} \omega_X^{[i]} = \mathbb{C}(Y)$ by multiplication by $\xi$. Suppose $f\in \mathbb{C}(Y)$ is integral over $\mathcal{O}_X$. For $0\leq i\leq m-1$, define $f_i = \frac{1}{m}\sum_{j=0}^{m-1} \xi^{-ij}g^if$. Then $f = \sum f_i$, where $f_i$ is integral and $f_i \in \mathbb{C}(X)\omega_X^{[i]}$. Then $f_i^m \in \omega_X^{[mi]}(mi\Delta)$, which implies $f_i \in \omega_X^{[i]}(\lfloor i\Delta \rfloor)$. 

(3) We follow the same argument as \cite[2.49]{Ko13}, which proved the local case. 

By (2), we have $p_*\mathcal{O}_Y = \oplus_{i=0}^{m-1} \omega_X^{[i]}(\lfloor i\Delta \rfloor)$. By duality for finite morphisms, see \cite[Proposition 5.68]{KM98}, we have $p_*\omega_Y = p_*\mathcal{H}om_Y(\mathcal{O}_Y, \omega_Y) = \mathcal{H}om_X(p_*\mathcal{O}_Y,\omega_X) = \oplus_{i=0}^{m-1} \omega_X^{[1-i]}(-\lfloor i\Delta \rfloor)$. Replacing $i$ by $1-i$, and observing that $\lfloor (1-i) \Delta\rfloor+ \lfloor i \Delta\rfloor = 0$ by Lemma \ref{Lem elementary}, we have $p_*\omega_Y  = \oplus_{i=-m+2}^{1} \omega_X^{[i]}(-\lfloor (1-i)\Delta \rfloor) =  \oplus_{i=-m+2}^{1} \omega_X^{[i]}(\lfloor i\Delta \rfloor) = \oplus_{i=0}^{m-1} \omega_X^{[i]}(\lfloor i\Delta \rfloor)$. 

So $p_*\omega_Y$ is isomorphic to $p_*\mathcal{O}_Y$, and this isomorphism preserves the $p_*\mathcal{O}_Y$-module structure. So $\omega_Y$ is isomorphic to $\mathcal{O}_Y$.

By \cite[Corollary 2.43]{Ko13}, $Y$ is klt. So $Y$ is canonical since $K_Y$ is Cartier.

(4) Let $\phi$ be an automorphism of $(X,\Delta)$, and $\phi^*$ be the corresponding automorphism of $\mathbb{C}(X)$. We have the following commutative diagram: 
\[\begin{tikzcd}
	Y & Y \\
	X & X
	\arrow["p", from=1-1, to=2-1]
	\arrow["p", from=1-2, to=2-2]
	\arrow["\phi", from=2-1, to=2-2]
    \arrow["\psi", dashed,  from=1-1, to=1-2]
\end{tikzcd}\]

Our goal is to construct the dashed arrow $\psi$. To achieve this, it suffices to show that the corresponding field extensions of $p$ and $\phi^{-1} \circ p$ are the same. 

Since $\phi$ is an automorphism of $(X,\Delta)$,  $\phi^*\omega_X \simeq \omega_X$, $\phi^*\Delta = \Delta$. So $\phi^*K_X \sim K_X$. Note that the index $1$ cover of $K_X+\Delta$ depends only on the linear equivalence class of $K_X$. But $p$ corresponds to the index $1$ cover of $K_X+\Delta$, while $\phi^{-1} \circ p$ correspond to the index $1$ cover of $\phi^*(K_X+\Delta) = \phi^*K_X + \Delta$. So our result holds. 

\hfill
\qedsymbol

\end{Pf}

\subsection{Polyhedral type cones}

Let $C$ be a cone, we denote the rational hull of $C$, i.e. the cone generated by rational points in the closure of $C$, by $C_+$. 

We recall the following result from \cite{Loo14}: 

\begin{Thm}\label{Prop-Def}
(\cite[Proposition-Definition 4.1]{Loo14}) 
Let $V$ be a real finite dimensional vector space
equipped with a rational structure $V(\mathbb{Q}) \subset V$ and let $C$ be an open nondegenerate convex cone in $V$. 

Let $\Gamma$ be a subgroup of $GL(V)$ which stabilizes
$C$ and some lattice in $V(\mathbb{Q})$. Then the following conditions are equivalent:

$(i)$ There exists a polyhedral cone $\Pi$ in $C_+$ with $\Gamma \cdot \Pi = C_+$.

$(ii)$ There exists a polyhedral cone $\Pi$ in $C_+$ with $\Gamma \cdot \Pi \supseteq C$.

Moreover, in case (ii) we necessarily have $\Gamma \cdot \Pi = C_+$. If one of these
equivalent conditions is fulfilled, we say that $(V(\mathbb{Q}), C, \Gamma)$ is a polyhedral
triple or simply, that $(C_+, \Gamma )$ is of polyhedral type.
\end{Thm}

\begin{Rem}
Moreover, when $(C_+, \Gamma )$ is of polyhedral type, by \cite[Application 4.14]{Loo14}, in fact there exists a rational polyhedral cone $\Pi \subset C_{+}$ which is a fundamental domain for the action of $\Gamma$ on $C_+$. 
\end{Rem}

\subsection{Variation of Hodge structures}

We briefly review some definitions and facts on variation of Hodge structures. 

We first recall the definition of Hodge structures: 

\begin{Def}(\cite[Definition 3.1.2]{CEGT14})
A Hodge structure of weight $k$ is defined by the following data:

1. A finitely generated abelian group $V_{\mathbb{Z}}$; 

2. A descending filtration $F^\bullet$ by complex subspaces $F^pV_{\mathbb{C}}$ of $V_{\mathbb{C}} = V_{\mathbb{Z}} \otimes_{\mathbb{Z}}\mathbb{C}$ satisfying
$$V_{\mathbb{C}} = F^pV_{\mathbb{C}}\oplus \overline{F^{k-p+1}}V_{\mathbb{C}}$$
We denote the $(p,q)$-part of $V_{\mathbb{C}}$ by $V_{\mathbb{C}}^{p,q} = F^pV_{\mathbb{C}}\cap \overline{F^qV_{\mathbb{C}}}$ for $p+q = k$. 

We denote a Hodge structure by $(V_{\mathbb{Z}}, F^\bullet)$. 

\end{Def}

\begin{Rem}
Sometimes $V_\mathbb{Z}$ is required to be free, which can be achieved easily by taking the quotient by the torsion subgroup. Our definition is sometimes called a $\mathbb{Z}$-Hodge structure in the literature. 
\end{Rem}

\begin{Def}(\cite[Definition 3.2.13]{CEGT14})\label{Def 3.2.13}
A polarization of a Hodge structure $(V_{\mathbb{Z}}, F^\bullet)$ of weight $k$ is a bilinear form $S : V_{\mathbb{Z}} \times V_{\mathbb{Z}} \to \mathbb{Z}$, whose extension by scalars is $S_{\mathbb{C}} : V_{\mathbb{C}} \times V_{\mathbb{C}} \to \mathbb{C}$, such that: 

(1) $S(x,y) = (-1)^k S(y,x)$ for $x,y \in V_{\mathbb{Z}}$ and $S(F^p,F^q) = 0$ for $p+q > k$, 

(2) (Hodge-Riemann relation)
$S_{\mathbb{C}}(C(V)u,\bar{v})$ is a positive definite Hermitian form on $V_{\mathbb{C}}$ where $C(V)$
 denotes the Weil action on $H$, i.e. $C(V)u := i^{p-q}u$ for $u \in V_{\mathbb{C}}^{p,q}$. 
    
\end{Def}

\begin{Rem}\label{Rem primitive}
It is well known that the cohomology $H^i(X,\mathbb{Z})$ of a nonsingular projective variety of dimension $d$ naturally carries a Hodge structure, see \cite[Theorem 3.1.14]{CEGT14}. 

If we fix an ample divisor $A$, then we have Lefschetz operator $L : H^i(X,\mathbb{Z}) \to H^{i+2}(X,\mathbb{Z})$.  The primitive part of $H^i(X,\mathbb{Z})$, defined as 
$$H^i_{prim}(X, \mathbb{Z}) = \ker(L^{d-i+1}: H^i(X,\mathbb{Z}) \to H^{2d-i+2}(X,\mathbb{Z}))$$
for $i\leq d$, 
carries a polarized Hodge structure. 

It is worth noting that if the middle cohomology $H^d(X,\mathbb{Z})$ is primitive, then we have a polarized Hodge structure on $H^d(X,\mathbb{Z})$ without giving an ample divisor $A$, since the polarization is given by the intersection pairing. 
\end{Rem}

Then we recall the definition of variation of Hodge structures, which describes a family of Hodge structures: 

\begin{Def}(\cite[Definition 7.3.4]{CEGT14}) \label{Def 7.3.4} 
Let $B$ be a connected complex manifold, a variation of Hodge structure of weight $k$ (VHS) over $B$ consists of a local system $\mathcal{V}_{\mathbb{Z}}$ of free $\mathbb{Z}$-modules and a filtration of the associated holomorphic vector bundle $\mathcal{V} = \mathcal{V}_{\mathbb{Z}}\otimes_{{\mathbb{Z}}} \mathcal{O}_B$:
$$\cdots \subset \mathcal{F}^p \subset \mathcal{F}^{p-1} \subset \cdots$$
by holomorphic subbundles $\mathcal{F}^p$ satisfying:

1. $\mathcal{V} = \mathcal{F}^p \oplus \overline{\mathcal{F}^{k-p+1}}$ as $C^\infty$ bundles, where the conjugation is relative to the local system of real vector spaces $\mathcal{V}_{\mathbb{R}} := \mathcal{V}_{\mathbb{Z}} \otimes \mathbb{R}$. 

2. Griffiths transversality: $ \nabla(\mathcal{F}^p) \subset \Omega^1_B \otimes  \mathcal{F}^{p-1}$, where $\nabla$ denotes the flat connection on $\mathcal{V}$, and we denote the vector bundle and the associated sheaf of sections by the same notation $\mathcal{F}^p$. 

We denote a variation of Hodge structure by $(\mathcal{V}_{\mathbb{Z}}, \mathcal{F}^{\bullet})$. 
\end{Def}

\begin{Rem}
For any $t\in B$, the restriction $(\mathcal{V}_{\mathbb{Z},t}, \mathcal{F}_t^{\bullet})$ is a Hodge structure. 
\end{Rem}

\begin{Def}
A polarization of a variation of Hodge structure $(\mathcal{V}_{\mathbb{Z}}, \mathcal{F}_{\bullet})$ over $B$ consists of a bilinear form $Q_{\mathcal{V}}: \mathcal{V}_{\mathbb{Z}} \times \mathcal{V}_{\mathbb{Z}} \to \mathbb{Z}$ such that: 

For any $t\in B$, the restriction of $Q_{\mathcal{V}}$ on $t$ gives a polarization of the Hodge structure $(\mathcal{V}_{\mathbb{Z},t}, \mathcal{F}_{t,\bullet})$. 
\end{Def}

It is known that a family of K\"ahler manifolds gives a variation of Hodge structure, and we restrict ourselves to algebraic case for simplicity: 

\begin{Thm}(\cite[10.2.1]{Voi02})\label{Thm geometric vhs}
Let $\phi: X \to B$ be a family of smooth projective varieties, i.e. $\phi$ is a smooth projective morphism between algebraic varieties. 

Then $\mathcal{V}_{\mathbb{Z}} : = R^k\phi_*\mathbb{Z}_X/\{Torsion\}$ is a local system. Let $\mathcal{V}$ be the associated holomorphic vector bundle $\mathcal{V}_{\mathbb{Z}}\otimes_{{\mathbb{Z}}} \mathcal{O}_B$. Then there exist holomorphic subbundles $\mathcal{F}^p$ such that 

(1)  $(\mathcal{F}^p)_b \subset (\mathcal{V})_b$ can be identified with $F^pH^k(X_b,\mathbb{C}) \subset  H^k(X_b,\mathbb{C})$ for every $ b \in B$, 

(2) $(\mathcal{V}_{\mathbb{Z}}, \mathcal{F}^\bullet)$ is a variation of Hodge structure. 
\end{Thm}

\begin{Rem}
This $(\mathcal{V}_{\mathbb{Z}}, \mathcal{F}^\bullet)$ is called the variation of Hodge structure associated to $\phi$. 
\end{Rem}

If we consider the primitive part of the cohomology (with respect to some polarization), we can further obtain a polarized variation of Hodge structure: 

\begin{Thm}(\cite[Section 2]{Gri70III})\label{Thm family gives polarized variation of Hs}
In the settings of Theorem \ref{Thm geometric vhs}, and let $\dim X/B = m$. Let $A$ be an ample divisor on $X$ over $B$, and $L: R^i\phi_*\mathbb{Z}_X \to R^{i+2}\phi_*\mathbb{Z}_X$ be the Lefschetz operator, which acts pointwise by $(-)\cup c_1(A|_{X_b}): H^k(X_b, \mathbb{Z}) \to H^{k+2}(X_b,\mathbb{Z})$. 

Let $0\leq k\leq m$, and $(R^k\phi_*\mathbb{Z}_X)_{prim}$ be the kernel of $L^{m-k+1}: R^{k}\phi_*\mathbb{Z}_X \to R^{2m-k+2}\phi_*\mathbb{Z}_X$. 

Then $\mathcal{V}_{\mathbb{Z}}^{prim} : = (R^k\phi_*\mathbb{Z}_X)_{prim}/\{Torsion\}$ is a local system. Let $\mathcal{V}^{prim}$ be the associated holomorphic vector bundle $\mathcal{V}_{\mathbb{Z}}^{prim}\otimes_{{\mathbb{Z}}} \mathcal{O}_B$. Then there exist holomorphic subbundles $\mathcal{F}^p$ such that 

(1) $(\mathcal{F}^p)_b \subset (\mathcal{V}^{prim})_b$ can be identified with $F^pH^k_{prim}(X_b,\mathbb{C}) \subset  H^k_{prim}(X_b,\mathbb{C})$ for every $ b \in B$, 

(2) $(\mathcal{V}_{\mathbb{Z}}^{prim}, \mathcal{F}^\bullet)$ is a variation of Hodge structure. 

(3) There exists a polarization $Q_{\mathcal{V}}$ of the variation of Hodge structure $(\mathcal{V}_{\mathbb{Z}}^{prim}, \mathcal{F}_{\bullet})$ such that for any $t\in B$, the restriction of $Q_{\mathcal{V}}$ on $t$ is given by $H^k_{prim}(X_b,\mathbb{C}) \times H^k_{prim}(X_b,\mathbb{C}) \to \mathbb{Z}$, $(u,v) \mapsto L^{m-k}(u \cup v)$. 
\end{Thm}

\begin{Rem}\label{Rem independence}
This $(\mathcal{V}_{\mathbb{Z}}^{prim}, \mathcal{F}^\bullet, Q_{\mathcal{V}})$ is called the polarized variation of Hodge structure associated to $(\phi,A)$. If moreover $k=m$ and $H^m(X_b,\mathbb{C})$ is primitive, $(\mathcal{V}_{\mathbb{Z}}^{prim}, \mathcal{F}^\bullet, Q_{\mathcal{V}})$ is independent of $A$ as in Remark \ref{Rem primitive}.  The ample divisor can be omitted in this special case, but the polarization pairing is
still part of the Hodge-theoretic data.
\end{Rem}

\subsection{Period domains and Period maps}

Following \cite[Section 8]{Gri70III}, we recall some basic facts on period domains and period maps. We also recommend \cite[Section 2,3]{Sch73} for a beautiful review. 

Let $E = E_\mathbb{R} \otimes_\mathbb{R} \mathbb{C}$ be a complex vector space with a real structure and 
$0<h_0\leq h_1\leq \cdots \leq h_{n-1}< h_n = \dim E$ 
an increasing sequence of integers such that 
$h_{n-q-1} = h_n - h_q$ for $0<q<n$. 
We also assume given a non-degenerate bilinear form 
$Q: E_\mathbb{R} \times E_\mathbb{R} \to \mathbb{R}$, $Q (e, e^\prime) = (- 1)^n Q(e^\prime, e)$, which extends $\mathbb{C}$-linearly to $E$. 

Consider the set $\check{D}$ of all filtrations 
$F^0 \subset \dots \subset F^n = E$, $\dim F^q = h_q$, 
which satisfy the first Riemann bilinear relations, i.e. $Q(F^q, F^{n-q-1}) = 0$. 

Consider the subset $\mathcal{D}$ of $\check{D}$, 
which satisfy the second Riemann bilinear relations, i.e. $Q(F^q, \overline{F^q})$ is non-degenerate, and $(-1)^q(-i)^nQ(E^q,\overline{E^q})$ is positive-definite. 

\begin{Prop}(\cite[Proposition 8.12]{Gri70III})
The period domain $\mathcal D$ is open in its compact dual $\check{D}$. Let
$G=\operatorname{Aut}(E_{\mathbb R},Q)$, and let $H$ be the stabilizer of a reference Hodge filtration. Then $H$ is compact and $\mathcal D\simeq G/H$.
\end{Prop}

\begin{Rem}
    The space $\mathcal{D}$ is called the period domain. 
\end{Rem}

Assume that we have a polarized variation of Hodge structure  $(\mathcal{V}_\mathbb{Z}, \mathcal{F}^\bullet, Q)$ on $B$, whose Hodge number is compatible with the period domain $\mathcal{D}$. Then we have a monodromy representation $\rho: \pi_1(B,b) \to Aut(\mathcal{V}_\mathbb{Z}|_b) \subset G$. The image of $\rho$ in $G$ is called the monodromy group of the polarized variation of Hodge structure $(\mathcal{V}_\mathbb{Z}, \mathcal{F}^\bullet, Q)$. 
We denote $\Gamma = \rho(G)$. 
Let $\Gamma_a := Aut(\mathcal{V}_\mathbb{Z}|_b,Q_b) \subset G$ be the arithmetic monodromy group, i.e. the automorphisms of the lattice $\mathcal{V}_\mathbb{Z}|_b$ that preserve the polarization $Q_b$. 

The quotient spaces $\Gamma \backslash \mathcal{D}$ and $\Gamma_a \backslash \mathcal{D}$ are separated analytic spaces by compactness of $H$. See \cite[(7.4.5)]{CEGT14} for more explanations. We have a quotient map $\Gamma \backslash \mathcal{D} \to \Gamma_a \backslash \mathcal{D}$. 

To state the correspondence between polarized variation of Hodge structures and period maps, we need some more definitions: 

\begin{Def}(\cite[Section 9(a)]{Gri70III})
A continuous mapping $\Phi : B\to \Gamma \backslash \mathcal{D}$ is said to be locally liftable if  for any $b \in B$, 
there exists an analytic neighborhood $U$ of $b$ and a continuous mapping $\tilde{\Phi} : U\to \mathcal{D}$ such that 
the following diagram is commutative: 
\[\begin{tikzcd}
	& {\mathcal{D}} \\
	U & {\Gamma \backslash\mathcal{D}}
	\arrow[from=1-2, to=2-2]
	\arrow["{\tilde{\Phi}}", from=2-1, to=1-2]
	\arrow["\Phi", from=2-1, to=2-2]
\end{tikzcd}\]
A locally liftable mapping $\Phi$ is holomorphic if the local liftings 
are holomorphic, and a locally liftable holomorphic mapping $\Phi$ is said to satisfy the 
infinitesimal period relation if the local liftings satisfy the infinitesimal period relation. 
\end{Def}

\begin{Rem}
The infinitesimal period relation is now called Griffiths transversality. 
This corresponds to Definition \ref{Def 7.3.4} (2) via the correspondence between variation of Hodge structures and period maps. 
We refer to \cite[(9.1)]{Gri70III} for precise formulation. 

\end{Rem}

Then we have the correspondence between polarized variation of Hodge structures and period maps. We only state one direction, which is enough for our purpose: 

\begin{Thm}(\cite[Proposition 9.3]{Gri70III})\label{Thm existence of period map}
A polarized variation of Hodge structure with monodromy 
group $\Gamma$ induces a locally liftable holomorphic mapping 
$\Phi: B \to \Gamma \backslash \mathcal{D}$
which satisfies the infinitesimal period relation. 
\end{Thm}

\begin{Rem}
We call the corresponding $\Phi$ the period map associated to the polarized variation of Hodge structure. Note that when we talk about the differentials of $\Phi$, we are really talking about the differentials of a local lifting. 
\end{Rem}

\subsection{Boundedness conventions}

\begin{Def}\label{Def boundedness conventions}
A collection $V$ of projective varieties is bounded if there is
a projective morphism $\mathcal X\to \mathcal S$, with $\mathcal S$ a scheme of finite
type over $\mathbb C$, such that every member in $V$ is isomorphic to a
closed fiber. A collection of fibrations $W$ is bounded if there is
a diagram of projective morphisms
$\mathcal X\to\mathcal B\to \mathcal S$, with $\mathcal S$ of finite type, such that
every fibration in $W$ is isomorphic to a fiber 
$\mathcal X_s\to\mathcal B_s$. For pairs, the isomorphisms must also
identify the boundaries with the restrictions of a fixed relative
divisor. 
\end{Def}

\begin{Def}\label{Def v-marked}
A $v$-marked curve fibration of a terminal $\mathbb{Q}$-factorial Calabi-Yau pair of dimension $d$ consists of the following data: 

(1) A terminal $\mathbb{Q}$-factorial Calabi-Yau pair $(X,\Delta)$ with a fibration $\pi: X\to B$ such that $\dim B = d - 1$; 

(2) A very ample divisor $A$ on $B$ such that $A^{d-1}\leq v$. 

We may write it as a 5-tuple $(X,\Delta,B,\pi,A)$. 

We say two $v$-marked curve fibrations of terminal $\mathbb{Q}$-factorial Calabi-Yau pairs $(X,\Delta,B,\pi,A)$ and $(X^\prime, \Delta^\prime, B^\prime, \pi^\prime, A^\prime)$ are isomorphic if there exists an isomorphism $(X, \Delta, B, \pi, [A]) \simeq (X^\prime, \Delta^\prime, B^\prime, \pi^\prime, [A^\prime])$, where $[-]$ denotes the numerical equivalence class. 

More precisely, an isomorphism $(X, \Delta, B, \pi, [A]) \simeq (X^\prime, \Delta^\prime, B^\prime, \pi^\prime, [A^\prime])$ consists of two isomorphisms $f: (X,\Delta) \to (X^\prime,\Delta^\prime)$, $g: B \to B^\prime$ such that $g\circ \pi = \pi^\prime \circ f$ and $g^*A^\prime$ is numerically equivalent to $A$. 
\end{Def}

\begin{Def}
A set $V$ of $v$-marked curve fibrations of terminal $\mathbb{Q}$-factorial Calabi-Yau pairs of dimension $d$ is said to be bounded if the following holds: 

There exist a pair $(\mathcal{X},\mathcal{D})$ and  quasi-projective varieties $\mathcal{B}, \mathcal{S}$ with a $\mathbb{Q}$-Cartier $\mathbb{Q}$-divisor $\mathcal{A}$ on $\mathcal{B}$
and a commutative diagram 
\[\begin{tikzcd}
	({\mathcal{X}}, {\mathcal{D}}) && {\mathcal{B}} \\
	& {\mathcal{S}}
	\arrow["f", from=1-1, to=1-3]
	\arrow["g"', from=1-1, to=2-2]
	\arrow["h", from=1-3, to=2-2]
\end{tikzcd}\]
of projective morphisms such that for any $v$-marked curve fibration of a terminal $\mathbb{Q}$-factorial Calabi-Yau pair $(X,\Delta,B,\pi,A) \in V$ there exists a closed point $s\in \mathcal{S}$ such that $(X,\Delta,B,\pi,A) \simeq (\mathcal{X}_s, \mathcal{D}_s, \mathcal{B}_s, f_s, [\mathcal{A}_s])$, where $[-]$ denotes the numerical equivalence class. 
\end{Def}

\begin{Rem}
We will only consider boundedness of the set of $v$-marked curve fibrations of strict terminal $\mathbb{Q}$-factorial Calabi-Yau pairs of dimension $3$. 
\end{Rem}

\section{Boundedness results for fibered Calabi-Yau pairs in dimension $3$}

In this section, we prove a generalization of boundedness of elliptic Calabi-Yau threefolds \cite[Theorem 1.1]{FHS21}. Our main result is Theorem \ref{Theorem Boundedness}. 

\subsection{Strict klt Calabi-Yau pairs}

\begin{Def}
Let $(X,\Delta)$ be a klt Calabi-Yau pair of dimension $3$ such that the Iitaka dimension $\kappa(-K_X) = 0$. 
We may run $K_X+(1+\epsilon)\Delta$-MMP on $X$ and terminate with a klt Calabi-Yau threefold $X^{min}$. 

We say that $(X,\Delta)$ is a strict klt Calabi-Yau pair if the augmented irregularity of $X^{min}$ is trivial. 

Furthermore, if $(X,\Delta)$ is terminal, we may call it a "strict terminal Calabi-Yau pair". If $\Delta = 0$, we call $X$ a strict klt Calabi-Yau variety. 

\end{Def}

\begin{Rem}
This definition is independent of the choice of $X^{min}$. 

Indeed, for another choice of MMP, say  $X \dashrightarrow \hat{X}^{min}$, $X^{min}$ and $\hat{X}^{min}$ are isomorphic in codimension $1$ since this MMP precisely contracts the support of $\Delta$. 

Moreover, for any quasi-\'etale cover $q: \tilde{X}^{min} \to X^{min}$, we have a corresponding quasi-\'etale cover $\hat{q}: \tilde{\hat{X}}^{min} \to \hat{X}^{min}$ obtained by taking the nomalization of $\hat{X}^{min}$ in the function field of $\tilde{X}^{min}$. In this case,  $\tilde{X}^{min}$ and $\tilde{\hat{X}}^{min}$ are isomorphic in codimension $1$. 
Since both $\tilde{X}^{min}$ and $\tilde{\hat{X}}^{min}$ have rational singularities, we have $h^1(\tilde{X}^{min},\mathcal{O}_{\tilde{X}^{min}}) = h^1(\tilde{\hat{X}}^{min},\mathcal{O}_{\tilde{\hat{X}}^{min}})$. 

\end{Rem}

\begin{Rem}
The assumption that $\kappa(-K_X) = 0$ implies that the coefficients of $\Delta$ lie in a uniform finite set, since the index of $K_{X^{min}}$ is bounded by \cite{Xu19b}. 

Indeed, by \cite[Theorem 1.13]{Xu19b}, there exists some integer $N$ such that $NK_{X^{min}} \sim 0$. Since $(X,\Delta) \dashrightarrow (X^{min},0)$ is crepant, we have $N(K_X+\Delta)\sim 0$, and in particular the coefficients of $\Delta$ lies in the finite set $\frac{1}{N}\mathbb{Z}\cap [0,1]$. 
\end{Rem}

\begin{Lem}\label{Lem strict klt CY pair}
Let $(X,\Delta)$ be a strict klt Calabi-Yau pair of dimension $3$. 

Then $h^1(X,\mathcal{O}_X) = h^2(X,\mathcal{O}_X) = 0$. 
\end{Lem}

\begin{Pf}
By definition, the Iitaka dimension $\kappa(-K_X) = 0$. We may run $K_X+(1+\epsilon)\Delta$-MMP on $X$ and terminate with a klt Calabi-Yau threefold $X^{min}$. Then by definition, the augmented irregularity of $X^{min}$ is trivial. 

Since $X$, $X^{min}$ have rational singularities and they are Cohen-Macaulay, see \cite[Corollary 2.88]{KM98}, it suffices to show $h^1(X^{min},\mathcal{O}_{X^{min}}) = h^2(X^{min},\mathcal{O}_{X^{min}}) = 0$. 

Let $\tilde{X}$ be the index $1$ cover of $X^{min}$, see Lemma \ref{Lemma delicate lifting property} with zero boundary. Since the augmented irregularity of $X^{min}$ is trivial, $h^1(\tilde{X},\mathcal{O}_{\tilde{X}})$ is trivial. By \cite[Proposition 5.75]{KM98}, the canonical sheaf $\omega_{\tilde{X}}$ coincides with the dualizing sheaf. 

By assumption, the canonical sheaf $\omega_{\tilde{X}}$ is isomorphic to $\mathcal{O}_{\tilde{X}}$. So Serre duality gives $h^2(\tilde{X},\mathcal{O}_{\tilde{X}}) = h^1(\tilde{X},\mathcal{O}_{\tilde{X}}) = 0$. Finally, $h^1(X^{min},\mathcal{O}_{X^{min}}) = h^2(X^{min},\mathcal{O}_{X^{min}}) = 0$ by \cite[Proposition 5.7 (2)]{KM98}. 

\hfill
\qedsymbol
\end{Pf}

\subsection{Boundedness of bases}

It is known that the set of K3 surfaces is not bounded, so the product of a K3 surface and an elliptic curve shows that the set of bases of elliptic Calabi-Yau pairs is not bounded. 

Nevertheless, we can prove the following result: 

\begin{Thm}\label{Thm boundedness of bases}
Let $V$ be the set of surfaces $B$ such that 

(1) There exists a klt Calabi-Yau pair $(X,\Delta)$ of dimension $3$, and a fibration $\pi: X\to B$, such that the generic fiber of $\pi$ is a curve;

(2) $h^1(X,\mathcal{O}_X) = h^2(X,\mathcal{O}_X) = 0$; 

(3) The coefficients of $\Delta$ lie in a finite set $\Phi \subset \mathbb{Q}$. 

Then $V$ is bounded. 

\end{Thm}

\begin{Pf}

By canonical bundle formula, we may write $K_X+\Delta \sim_{\mathbb{Q}} f^*(K_B+\Delta_B+M_B)$, where $(B,\Delta_B+M_B)$ is a generalised klt pair, see \cite[Definition 4.1]{BZ16} for definition of generalised pairs. By \cite[Remark 5.1]{Fil20}, we may choose $0 \leq \Delta_B^\prime \sim_{\mathbb{Q}} M_B$ such that $(B, \Delta_B +\Delta_B^\prime)$ is klt, and the coefficients of $\Delta_B +\Delta_B^\prime$ belong to a DCC set of rational numbers only depending on $n$ and $\Phi$. By global ACC, the coefficients of $\Delta_B +\Delta_B^\prime$ lie in a finite set by \cite[Theorem 1.5]{HMX12}. By \cite[Theorem 1.5]{Xu19a}, the index of $K_X+\Delta_B +\Delta_B^\prime$ is bounded, in particular, $X$ is $\epsilon$-lc for some fixed $\epsilon$. 

We claim that $h^i(B,\mathcal{O}_B) = 0$ for $i=1,2$. By Hodge symmetry for klt pairs, see \cite[Theorem 1]{Sch16}, $h^i(B,\mathcal{O}_B) = h^0(B,\Omega_B^{[i]})$ for $i=1,2$. By \cite[Theorem 2]{Sch16}, we may pull back the differential forms to $X$, so we have an injection $H^0(B,\Omega_B^{[i]}) \to H^0(X,\Omega_X^{[i]}) = H^i(X,\mathcal{O}_X) = 0$ for $i=1,2$. So our claim holds. 

If $\Delta_B+M_B$ is nonzero or $B$ is not canonical, then $B$ is uniruled, 
hence is rational or ruled. When $B$ is rational, $B$ is bounded by \cite[Theorem 1.6]{Bir23}. When $B$ is birational to a non-rational ruled surface, 
$h^1(B,\mathcal{O}_B) \neq 0$, which cannot appear. 

Otherwise, we have $K_B \sim_{\mathbb{Q}} 0$, and $B$ is canonical. Let $S$ be the minimal resolution of $B$. Since $B$ has rational singularities, $h^i(S,\mathcal{O}_S) = h^i(B,\mathcal{O}_B) = 0$ for $i=1,2$. By classification of surfaces, $S$ must be an Enriques surface. Then all possible $B$ are bounded by \cite[Theorem 6.3]{FHS21}. 

\hfill
\qedsymbol
\end{Pf}

\subsection{Boundedness of rational multi-sections}

In order to prove boundedness of curve fibered Calabi-Yau pairs in dimension $3$, we also need boundedness of rational multi-sections.

For a smooth genus-zero generic fiber, there is a direct uniform bound:

\begin{Prop}\label{Proposition multi-section 0}
Let $\pi:X\to B$ be a fibration whose generic fiber is a smooth curve of genus zero. Then $\pi$ admits a
rational multisection of degree at most two.
\end{Prop}

\begin{Pf}
Put $K=\mathbb C(B)$. The generic fiber is a smooth conic over $K$.
Its anticanonical linear system has an effective divisor of degree
two defined over $K$, so it contains a closed point of degree one or
two. The closure of that point in $X$ dominates $B$ generically
finitely with the same degree, and gives the required rational
multisection.
\hfill\qedsymbol
\end{Pf}

Here an \emph{elliptic fibration} means a fibration whose generic
fiber is a smooth geometrically integral curve of genus one; a
section is not required. We now consider such fibrations: 

\begin{Prop}\label{Proposition multi-section 1}
Let $X$ be a strict klt Calabi-Yau threefold, and let $\pi: X \to B$ be an elliptic fibration. 

Then there exists an integer $d$  such that $\pi$ admits a rational multi-section of degree $d$. 
\end{Prop}

\begin{Pf}
Let $\mu: \tilde{X} \to X$ be the index $1$ cover of $X$, and $p: Y\to \tilde{X}$ a $\mathbb{Q}$-factorial terminalization. Then the degree of $\mu$ is bounded by \cite[Theorem 2.16]{Xu19b}. 

Let $Y \to C \to B$ be the Stein factorization of $\pi\circ \mu\circ p$. Then $Y\to C$ is an elliptic fibration of terminal $\mathbb{Q}$-factorial Calabi-Yau threefold in the sense of \cite{FHS21}, so it lies in a bounded family by \cite[Theorem 1.1]{FHS21}. In particular, $Y\to C$ admits a rational multi-section of bounded degree. Moreover, the degree of $C\to B$ is bounded by the degree of $\mu$, so $Y\to B$ admits a rational multi-section $\sigma$ of bounded degree. 

Then $\mu \circ p\circ \sigma$ exhibits a multi-section of $\pi$ of bounded degree. 

\hfill
\qedsymbol
\end{Pf}

\begin{Prop}\label{Proposition multi-section 2}
Let $(X,\Delta)$ be a strict klt Calabi-Yau pair of dimension $3$, and let $\pi: X \to B$ be an elliptic fibration. 

Then there exists an integer $d$  such that $\pi$ admits a rational multi-section of degree $d$. 
\end{Prop}

\begin{Pf}
Let $(X,\Delta) \to B$ be an elliptic fibration. By assumption, $\Delta$ is vertical, we may run $K_X+(1+\epsilon) \Delta$-MMP on $X$ over $B$, which terminates with a good minimal model $(Y,(1+\epsilon)\Delta_Y)$, inducing an Iitaka fibration $f: Y \to C$. In sum, we have: 

\[\begin{tikzcd}
	{(X,\Delta)} & {(Y,\Delta_Y)} \\
	B & C
	\arrow[dashed, from=1-1, to=1-2]
	\arrow[from=1-1, to=2-1]
	\arrow[from=1-2, to=2-2]
	\arrow[from=2-2, to=2-1]
\end{tikzcd}\]

Since $\Delta_Y$ is vertical, $g: C\to B$ is birational. By definition of Iitaka fibration, there exists some $\mathbb{Q}$-Cartier $\mathbb{Q}$-divisor $\Delta_C$ on $C$ such that $\Delta_Y \sim_{\mathbb{Q},B} f^*\Delta_C$. We may further assume $\Delta_Y \sim_{\mathbb{Q}} f^*\Delta_C$ by adding some $\mathbb{Q}$-Cartier $\mathbb{Q}$-divisor of the form $g^*D$, where $D$ is a $\mathbb{Q}$-Cartier $\mathbb{Q}$-divisor on $B$. Since $\Delta\geq0$ and $\kappa(\Delta)=0$, one has $h^0(X,\mathcal O_X(m\Delta))=1$ for sufficiently divisible $m>0$, as two independent sections would yield a nonconstant rational map and hence positive Iitaka dimension. The above MMP is a $\Delta$-MMP, so its steps preserve these section spaces. Indeed, on a common resolution of a step, the pullback of the initial divisor equals the pullback of its transform plus an effective exceptional divisor; pushing forward the corresponding sheaves gives the equality of sections. Since $f_*\mathcal O_Y=\mathcal O_C$, the projection formula and $\Delta_Y\sim_{\mathbb Q}f^*\Delta_C$ give \[ H^0(X,\mathcal O_X(m\Delta))  \simeq H^0(Y,\mathcal O_Y(m\Delta_Y))  \simeq H^0(C,\mathcal O_C(m\Delta_C)). \] All three spaces therefore have dimension one for sufficiently divisible $m$. Thus $\Delta_C$ has a unique effective representative up to $\mathbb Q$-linear equivalence. Replacing it by that representative, we may assume $\Delta_C\geq0$. 

Now we claim that $\Delta_Y = f^*\Delta_C$. Indeed, both of them are effective $\mathbb{Q}$-Cartier $\mathbb{Q}$-divisors, and they are known to be $\mathbb{Q}$-linearly equivalent, so they must be the same as $h^0(Y,\mathcal{O}_Y(m\Delta_Y)) = 1$ for sufficiently divisible $m$. 
Moreover, we have $\kappa(\Delta_C) = \kappa(\Delta_Y) = \kappa(\Delta) = 0$, which still follows from the fact that $h^0(X,\mathcal{O}_X(m\Delta)) = h^0(Y,\mathcal{O}_Y(m\Delta_Y)) = h^0(C,\mathcal{O}_C(m\Delta_C)) = 1$ for sufficiently divisible $m$. 

Note that $C$ admits a boundary divisor $D_C$ such that $(C,D_C)$ is a klt Calabi-Yau pair by \cite[Remark 5.1]{Fil20}. We may run $K_C+D_C+\epsilon\Delta_C$-MMP, which terminates with $C^{min}$, such that $C\to C^{min}$ contracts all the components of $\Delta_C$. Since $\Delta_Y$ is a very exceptional divisor over $C^{min}$ in the sense of \cite[Definition 3.1]{Bir11}, by \cite[Theorem 1.8]{Bir11}, we may run $K_Y+(1+\epsilon) \Delta_Y$-MMP on $Y$ over $C^{min}$, and terminate with $Y^{min}$ such that $K_Y^{min} \sim_{\mathbb{Q}} 0$, and $\Delta_{Y^{min}} = 0$. In sum, we have: 
\[\begin{tikzcd}
	{(X,\Delta)} & {(Y,\Delta_Y)} & {(Y^{min},0)} \\
	B & C & {C^{min}}
	\arrow[dashed, from=1-1, to=1-2]
	\arrow[from=1-1, to=2-1]
	\arrow[dashed, from=1-2, to=1-3]
	\arrow[from=1-2, to=2-2]
	\arrow[from=1-3, to=2-3]
	\arrow[from=2-2, to=2-1]
	\arrow[from=2-2, to=2-3]
\end{tikzcd}\]

Since horizontal arrows are birational, and rational multi-sections are preserved by birational equivalences, we may replace $(X,\Delta) \to B$ by $(Y^{min},0) \to C^{min}$. Then the result follows from Proposition \ref{Proposition multi-section 1}. 

\hfill
\qedsymbol
\end{Pf}

\subsection{Boundedness of marked curve fibered Calabi-Yau pairs}

In order to deduce finiteness of curve fibrations up to automorphism, we need boundedness of curve fibered Calabi-Yau pairs. And then one might apply some discreteness arguments. 

For surface fibrations, we are not able to prove boundedness of surface fibered Calabi-Yau pairs.
Nevertheless, it turns out that a stronger version of boundedness of curve fibered Calabi-Yau pairs leads to finiteness of surface fibrations up to automorphism, see Definition \ref{Def v-marked}.

We need the following lemma, which is a slight modification of \cite[Theorem 4.2]{FHS21}, and the proof is exactly the same: 

\begin{Lem}{\cite[Theorem 4.2]{FHS21}}\label{Lem deformation of divisor}
Let $(X_0,\Delta_0)$ be a terminal $\mathbb{Q}$-factorial pair. 
Assume that $K_{X_0} + \Delta_0 \equiv 0$, $H^1(X_0, \mathcal{O}_{X_0}) = 0$ and $H^2(X_0, \mathcal{O}_{X_0}) = 0$. 
Given a deformation $(X,\Delta) \to (T, 0)$ of $(X_0,\Delta_0)$ over a smooth variety $T$ , then $(X,\Delta)$ is terminal
$\mathbb{Q}$-factorial, $K_X+\Delta \sim_{\mathbb{Q},T} 0$ over a neighborhood of $0 \in T $. 
Furthermore, after an \'etale base change, the following
facts hold:

(1) $\overline{A}(X_\eta) \supset \overline{A}(X_0)$;

(2) $\overline{B}(X/T ) \subset \overline{B}(X_0)$; and

(3) $\overline{M}(X/T ) \supset \overline{M}(X_0)$.

\end{Lem}

Then we prove boundedness of $v$-marked curve fibrations of strict terminal $\mathbb{Q}$-factorial Calabi-Yau pairs of dimension $3$: 

\begin{Thm}\label{Theorem Boundedness}
Fix a positive real number $v$. The set of $v$-marked curve fibrations of strict terminal $\mathbb{Q}$-factorial Calabi-Yau pairs of dimension $3$ forms a bounded family. 
\end{Thm}

\begin{Pf}
\noindent \textbf{Step 1: boundedness in codimension 1}

By Proposition \ref{Proposition multi-section 0} and Proposition \ref{Proposition multi-section 2}, there exists a positive integer $C$, independent of the 5-tuple $(X, \Delta, B, \pi, A)$, such that $\pi$ admits a rational $d$-section, for
some $d \leq C$. By the remark after Lemma \ref{Lem strict klt CY pair}, the coefficients of $\Delta$ lie in a finite set $\Phi$. Possibly enlarging $\Phi$, we may assume $0,\frac{1}{2} \in \Phi$. 

Let $(X,\Delta,B,\pi,A)$ be a $v$-marked curve fibration of strict terminal $\mathbb{Q}$-factorial Calabi-Yau pair of dimension $3$. Then we may replace $A$ by a general member in its linear system. Now consider the morphism of pairs $\pi: (X,\Delta+\frac{1}{2}\pi^*A) \to (B, \Delta_B)$, where $\Delta_B$ is given by taking canonical bundle formula with respect to $\pi: (X,\Delta+\frac{1}{2}\pi^*A) \to B$ and then choosing a suitable representative of the moduli part as in \cite[Remark 5.1]{Fil20}. 

In the notation of \cite[Theorem 7.2]{Fil20}, $\pi: (X,\Delta+\frac{1}{2}\pi^*A) \to (B, \Delta_B)$ is an element in $\mathfrak{C}(3, v, \Phi, C)$, so the set of such $\pi$ forms a bounded family in codimension $1$, in the sense of \cite[Definition 2.6]{Fil20}: 

There exist a pair $(\mathcal{X},\mathcal{D})$ and quasi-projective varieties $\mathcal{B}$, $\mathcal{S}$ with a reduced Weil divisor $\mathcal{A}$ on $\mathcal{X}$ and a commutative diagram
\[\begin{tikzcd}
	{\mathcal{X}} && {\mathcal{B}} \\
	& {\mathcal{S}}
	\arrow["f", from=1-1, to=1-3]
	\arrow["g"', from=1-1, to=2-2]
	\arrow["h", from=1-3, to=2-2]
\end{tikzcd}\]
of projective morphisms such that for any morphism $\pi: (X,\Delta+\frac{1}{2}\pi^*A) \to B$ constructed as above, there exists some closed point $s\in \mathcal{S}$ such that there exists a commutative diagram
\[\begin{tikzcd}
	{\mathcal{X}_s} & X \\
	{\mathcal{B}_s} & B
	\arrow["{\mu_s}", dashed, from=1-1, to=1-2]
	\arrow["{f_s}"', from=1-1, to=2-1]
	\arrow["\pi", from=1-2, to=2-2]
	\arrow["{\nu_s}", from=2-1, to=2-2]
\end{tikzcd}\]
where $\nu_s$ is an isomorphism, $\mu_s$ is an isomorphism in codimension $1$, and $\text{Supp}(\mu_s^{-1}A) \subset \text{Supp}(\mathcal{A}_s)$, $\text{Supp}(\mu_s^{-1}\Delta) \subset \text{Supp}(\mathcal{D}_s)$. From the proof of \cite[Theorem 7.2]{Fil20}, we may further assume the $\mathcal{X}_s$ above is $\mathbb{Q}$-factorial.  
Note that the boundary divisor on the base is ignored here. 

\medskip
\noindent \textbf{Step 2: descent of the boundary}

By passing to a flattening stratification, we may assume $g$ and $h$ are flat. As in \cite[Proposition 2.10]{BDCS20}, by passing to a stratification and possibly discarding some components, we may assume that:

1) $\mathcal{S}$ is smooth,

2) $\mathcal{A}$ does not contain any fibre of $g$, 

3) every fibre $\mathcal{X}_s$ is a normal variety, and

4) for any $s \in \mathcal{S}$ and any irreducible component $\mathcal{A}^\prime$ of $\mathcal{A}$, $\mathcal{A}^\prime_s$ is either empty or an irreducible prime divisor on $\mathcal{X}_s$.

By copying the fibrations $\mathcal{X}\to \mathcal{B}\to \mathcal{S}$ finitely many times and equipping each with one single irreducible component of $\mathcal{A}$ and a subset of irreducible components of $\mathcal{D}$ with coefficients in $\Phi$, we may assume further that the above inclusion $\text{Supp}(\mu_s^{-1}A) \subset \text{Supp}(\mathcal{A}_s)$ is an equality, the coefficient of $\mathcal{A}$ is $1/2$ and $\mu_s^{-1}\Delta = \mathcal{D}_s$.

After an \'etale base change and a stratification, we may assume Lemma \ref{Lem deformation of divisor} holds, i.e. for each connected component $U$ of $\mathcal{S}$, there exists a distinguished point $0$ such that Lemma \ref{Lem deformation of divisor} holds for $(U,0)$.

Up to a stratification, we may also assume the set of points $s\in \mathcal{S}$ such that the fiber over $s$ is isomorphic in codimension $1$ to some $v$-marked curve fibrations of strict terminal $\mathbb{Q}$-factorial Calabi-Yau pair is dense. 

Since $\mathcal{X}$ is $\mathbb{Q}$-factorial, $\mathcal{A}$ is $\mathbb{Q}$-Cartier. For any ample divisor $\mathcal{H}$ on $\mathcal{X}$ and any $n\in \mathbb{Z}$, $\mathcal{H}+n\mathcal{A}$ is ample over $\mathcal{B}$ since this holds for a dense set of fibers and ampleness is an open condition. This means $\mathcal{A}$ is numerically trivial over $\mathcal{B}$. By the same argument, $K_\mathcal{X}+\mathcal{D}$ is also numerically trivial over $\mathcal{B}$. 

Let $c = 0,1$. Since $(\mathcal{X}_s, \mathcal{D}_s+\frac{c}{2}\mathcal{A}_s)$ is klt for a dense subset of $\mathcal{S}$, by iterated inversion of adjunction, $(\mathcal{X},\mathcal{D}+\frac{c}{2}\mathcal{A})$ is klt. By \cite[Corollary 1.6]{HX11}, $K_\mathcal{X}+\mathcal{D}+\frac{c}{2}\mathcal{A}$ is $\mathbb{Q}$-linearly trivial over $\mathcal{B}$. Taking the difference for $c = 0,1$, we deduce that $\mathcal{A}$ is $\mathbb{Q}$-linearly trivial over $\mathcal{B}$. 

In particular, there exists an ample $\mathbb{Q}$-divisor $\mathcal{H}$ on $\mathcal{B}$ such that $f^*\mathcal{H} = \mathcal{A}$. 

\medskip
\noindent \textbf{Step 3: application of the cone conjecture}

In this step, we essentially follow the same argument as \cite[Theorem 6.18]{FHS21}. 

For any $v$-marked curve fibration of
a strict terminal $\mathbb{Q}$-factorial Calabi-Yau pair of dimension $3$ $(X,\Delta,B,\pi,A)$, there exists an irreducible component $\mathcal{T}$ of $\mathcal{S}$ such that for some closed point $t\in \mathcal{T}$, $(X^\prime,\Delta^\prime,B,\pi^\prime,[A]) \simeq (\mathcal{X}_t, \mathcal{D}, \mathcal{B}_t, f_t, [\mathcal{H}_t])$, where $(X^\prime,\Delta^\prime,\pi^\prime)$ is isomorphic in codimension $1$ with $(X,\Delta,\pi)$ over $B$. 

By \cite[Theorem 1.4]{Li23}, $(\mathcal{X}_{\mathcal{T}},\mathcal{D}_{\mathcal{T}}) \to \mathcal{B}_{\mathcal{T}}$ admits only finitely many minimal models $(\mathcal{X}_{\mathcal{T},1},\mathcal{D}_{\mathcal{T},1})$, \dots , $(\mathcal{X}_{\mathcal{T},k},\mathcal{D}_{\mathcal{T},k})$ over $\mathcal{B}_{\mathcal{T}}$, up to isomorphism over $\mathcal{B}_{\mathcal{T}}$. Since $X$ and $X^\prime$ are $\mathbb{Q}$-factorial and isomorphic in codimension $1$, there exists some divisor $H$ which lies in the closure of the movable cone of $X^\prime/B$ such that the MMP of $(X^\prime,\Delta+\epsilon H)$ over $B$ terminates with $X$, where $\epsilon$ is a small enough positive number. This $H$ can be chosen to be the strict transform of some ample divisor $C$ on $X$. We may view $H$ as a divisor on $\mathcal{X}_t$ and extend to a divisor $\mathcal{H}$ on $\mathcal{X}_{\mathcal{T}}$ which lies in the closure of the movable cone of $\mathcal{X}_{\mathcal{T}}/ \mathcal{B}_{\mathcal{T}}$ by Lemma \ref{Lem deformation of divisor}. 

Then we run MMP of $(\mathcal{X}_{\mathcal{T}}, \mathcal{D}+\epsilon \mathcal{H})$ over $\mathcal{B}_{\mathcal{T}}$, which terminates with one of $(\mathcal{X}_{\mathcal{T},i},\mathcal{D}_{\mathcal{T},i})$, $1\leq i\leq k$. By \cite[Section 3]{HMX18a}, this MMP restricts to a sequence of $H$-negative birational maps whose composite an isomorphism in codimension $1$. So the induced birational map $\phi: (X,\Delta) \dashrightarrow (\mathcal{X}_{\mathcal{T},i,t}, \mathcal{D}_{\mathcal{T},i,t})$ is an isomorphism in codimension $1$, and the strict transform of the ample divisor $C$ on $\mathcal{X}_{\mathcal{T},i,t}$ is nef. By negativity lemma, this implies $\phi$ is crepant and $\phi^{-1}$ is a morphism. But then $\phi$ must be an isomorphism since $X$ is $\mathbb{Q}$-factorial. 

Hence we can give a family bounding all $v$-marked curve fibrations of
strict terminal $\mathbb{Q}$-factorial Calabi-Yau pairs of dimension $3$ by taking all relative minimal models of all components of $\mathcal{X} \to \mathcal{B} \to \mathcal{S}$. 

\hfill
\qedsymbol

\end{Pf}

\begin{Rem}
By Theorem \ref{Thm boundedness of bases}, the set of bases of curve fibrations of strict terminal $\mathbb{Q}$-factorial Calabi-Yau pairs forms a bounded family, so there exists a universal constant $v_0$ such that any curve fibration of a strict terminal $\mathbb{Q}$-factorial Calabi-Yau pair of dimension $3$ admits a $v_0$-marking. 
Forgetting the $v_0$-marking in Theorem
\ref{Theorem Boundedness} recovers boundedness of the underlying curve fibrations of strict terminal $\mathbb Q$-factorial
Calabi–Yau threefold pairs. 
\end{Rem}

\section{Period maps and isomorphism loci of strict Calabi--Yau pairs}

In this section, we prove an extraction result of strict terminal $\mathbb{Q}$-factorial Calabi-Yau pairs by analysis of period maps. Our main result is Theorem \ref{Thm extraction of strict Calabi-Yau pairs}. 

\subsection{Some Hodge theoretic properties of strict terminal Calabi-Yau pairs}

Let $X$ be a terminal complex threefold. Write $IH^j(X,\mathbb Q)$ for middle-perversity intersection cohomology. Then we have a comparison map $H^3(X,\mathbb{Q}) \to IH^3(X,\mathbb{Q})$ induced by the natural morphism from the constant sheaf to the intersection complex. 

\begin{Lem}\label{Lem analytic Q-factorization}
Let $X$ be a terminal variety of dimension $3$, then the natural map $H^3(X,\mathbb{Q}) \to IH^3(X,\mathbb{Q})$ is surjective. 
\end{Lem}

\begin{Pf}
Let $\pi:Y\to X$ be a proper small analytic $\mathbb Q$-factorialization, which exists by \cite[Corollary 4.5$^\prime$]{Kaw88}. By \cite[Corollary 4.12]{Kol89}, $\pi$ induces an isomorphism $IH^3(X,\mathbb{Q}) \to IH^3(Y,\mathbb{Q})$. 

Recall that a complex threefold $Z$ is a rational homology manifold if, for every $z\in Z$, its local rational homology
$H_i(Z,Z\backslash \{z\};\mathbb Q)$ is $\mathbb Q$ for $i=6$ and zero otherwise. 
By \cite[Lemma 4.2]{Kol89}, $Y$ has this property. The hypothesis in \cite[Lemma 4.2]{Kol89} concerns the analytic local rings, so an algebraic $\mathbb Q$-factorialization from \cite{BCHM} alone would not enough for our purpose.
The homology-manifold property gives $H^3(Y,\mathbb Q)\simeq IH^3(Y,\mathbb Q)$. 

Consider the following commutative diagram: 

\[\begin{tikzcd}
	{H^3(Y,\mathbb{Q})} & {IH^3(Y,\mathbb{Q})} \\
	{H^3(X,\mathbb{Q})} & {IH^3(X,\mathbb{Q})}
	\arrow["\simeq", from=1-1, to=1-2]
	\arrow[from=2-1, to=1-1]
	\arrow[from=2-1, to=2-2]
	\arrow["\simeq"', from=2-2, to=1-2]
\end{tikzcd}\]

It suffices to show the natural map $H^3(X,\mathbb{Q}) \to H^3(Y,\mathbb{Q})$ is surjective. 

Consider the spectral sequence $H^i(X,R^j\pi_*\mathbb{Q}_Y) \Rightarrow H^{i+j}(Y,\mathbb{Q}_Y)$. Since $\pi$ is small, it follows that $H^i(X,R^j\pi_*\mathbb{Q}_Y) = 0$ for $i>0,j>0$ or $j>2$, so $H^0(X,R^3\pi_*\mathbb{Q}_Y) = H^1(X,R^2\pi_*\mathbb{Q}_Y) = H^2(X,R^1\pi_*\mathbb{Q}_Y) = 0$. So the natural map $H^3(X,\mathbb{Q}) \to H^3(Y,\mathbb{Q})$, which is identified with the edge morphism of the $E_2$ page of the spectral sequence, is surjective. 

\hfill
\qedsymbol
\end{Pf}

\begin{Lem}\label{Lem cohomology of Calabi-Yau type}
Let $X$ be a terminal variety of dimension $3$ such that $h^0(X,\omega_X) = 1$, and $\mu: Y\to X$ be a resolution of singularities. Then the image of the natural map $H^3(X,\mathbb{C}) \to H^3(Y,\mathbb{C})$ contains $H^{3,0}(Y,\mathbb{C})$. 
\end{Lem}

\begin{Pf}
Since $X$ has rational singularities by \cite[Corollary 2.88]{Ko13},  $h^0(Y,\omega_Y) = 1$ by \cite[Definition 2.76 (3)]{Ko13}. 

Note that $H^3(X,\mathbb{Q}) \to H^3(Y,\mathbb{Q})$ is a morphism of mixed Hodge structures, and it is injective in weight $3$ part by \cite[Theorem 11.1.19 (iv)]{Max19}. Since $H^{3,0}(Y,\mathbb{C})$ is of dimension $1$, it suffices to show that $(3,0)$-part of $Gr_3^WH^3(X,\mathbb{C})$ is nonzero. By Lemma \ref{Lem analytic Q-factorization}, $H^3(X,\mathbb{C})$ maps surjectively to $IH^3(X,\mathbb{C})$. So it suffices to show $(3,0)$-part of $IH^3(X,\mathbb{C})$ is nonzero. 

Let $U$ be the smooth locus of $X$. Possibly passing to a higher resolution, we may assume the exceptional locus is  snc, say $E$. By \cite[Theorem 11.4.4]{Max19}, we have $IH^3(X,\mathbb{C}) \simeq Gr_3^WH^3_c(U,\mathbb{C})$. 
Identify $U=X_{\mathrm{reg}}$ with its inverse image in $Y$. Consider the exact sequence of mixed Hodge structures associated with $U=Y\backslash  E$: 
\[
\cdots\longrightarrow H_c^3(U,\mathbb C)
\longrightarrow H^3(Y,\mathbb C)
\longrightarrow H^3(E,\mathbb C)\longrightarrow\cdots.
\]
By \cite[Corollary 3.2.33]{CEGT14}, $Gr^W_3H^3(E)$ does not contain $(3,0)$-part, since it is a subquotient of Hodge structures given by varieties of dimension at most $2$. So we conclude our result. 

\hfill
\qedsymbol
\end{Pf}

\begin{Lem}\label{Lem intersection form}
Let $X$ be a terminal variety of dimension $3$, and $\mu: Y\to X$ be a resolution of singularities. Let $V$ be the image of the natural map $H^3(X,\mathbb{C}) \to H^3(Y,\mathbb{C})$. 

Then the intersection form of $H^3(Y,\mathbb{C})$ is nondegenerate when restricted to $V$, and this defines a bilinear form on $V$, which is independent of the choice of the resolution. 
\end{Lem}

\begin{Pf}
Let $\pi: Z \to X$ be an analytic $\mathbb{Q}$-factorization, which exists by \cite[Corollary 4.5$^\prime$]{Kaw88}, and $W$ be a common resolution of $Z$ and $Y$, with natural maps $p:W\to Y$, $q:W\to Z$. We note that by \cite[Theorem 5]{Moi67}, $W$ can be chosen to be projective, while $Z$ might not be projective. 

Then we have a commutative diagram: 
\[\begin{tikzcd}
	& {H^3(W,\mathbb{C})} \\
	{H^3(Y,\mathbb{C})} && {H^3(Z,\mathbb{C})} \\
	& {H^3(X,\mathbb{C})}
	\arrow["{p^*}", from=2-1, to=1-2]
	\arrow["{q^*}"', from=2-3, to=1-2]
	\arrow["{\mu^*}", from=3-2, to=2-1]
	\arrow["{\pi^*}"', from=3-2, to=2-3]
\end{tikzcd}\]

Then $H^3(Y,\mathbb{C})$, $H^3(Z,\mathbb{C})$, $H^3(W,\mathbb{C})$ admit natural intersection forms, and $p^*$, $q^*$ are compatible with them. 

So it suffices to show the intersection form is nondegenerate when restricted to $p^*V = q^*\pi^*H^3(X,\mathbb{C}) = q^*H^3(Z,\mathbb{C})$, since $\pi^*$ is surjective by Lemma \ref{Lem analytic Q-factorization}. But this follows immediately from the compatibility of $q^*$ with the intersection form and the fact that Poincar\'e duality holds for the $\mathbb{Q}$-homology manifold $Z$. 

To show independence of resolution, let $\mu_1: Y_1\to X$, $\mu_2:Y_2\to X$ be two resolutions, and $Y_0$ be a common resolution, with $p_1: Y_0\to Y_1$, $p_2:Y_0\to Y_2$. Then the intersection forms defined by $\mu_1$ and $\mu_2$ are the same by compatibility of $p_1^*$, $p_2^*$ with the intersection form. 

\hfill
\qedsymbol
\end{Pf}

\begin{Lem}\label{Lemma independence of resolution}
Let $X$ be a projective terminal variety of dimension $3$ with $H^1(X,\mathcal O_X)=0$, and let $\mu:Y\to X$ be a  resolution of singularities. Put
\[
V_{\mathbb Z}=
\operatorname{Im}\bigl(H^3(X,\mathbb Z)\longrightarrow
H^3(Y,\mathbb Z)\bigr)/\mathrm{torsion},
\qquad V_{\mathbb C}=V_{\mathbb Z}\otimes\mathbb C.
\]
The induced filtration $F_V^p=F^pH^3(Y,\mathbb C)\cap V_{\mathbb C}$ and the restriction of the intersection pairing give a polarized Hodge structure on $V_{\mathbb Z}$, which is independent of the resolution.
\end{Lem}

\begin{Pf}
Terminal singularities are rational, so
$H^1(Y,\mathcal O_Y)=H^1(X,\mathcal O_X)=0$.
Consequently $H^1(Y,\mathbb C)=0$, so $H^3(Y,\mathbb Q)$ is primitive. Its intersection pairing is therefore a polarization in the convention of Definition \ref{Def 3.2.13}.

The pullback on cohomology is a morphism of mixed Hodge structures. Its image $V_{\mathbb Q}$ is thus a pure Hodge substructure of $H^3(Y,\mathbb Q)$. Restriction of the associated positive definite Hermitian form to this substructure remains positive definite; hence the restricted bilinear form gives a polarization.

Then we prove independence of the resolution. The Hodge structure is independent of the  resolution because this filtration coincides with the Hodge structure of the weight $3$ part of the mixed Hodge structure of $H^3(X,\mathbb{Z})$ by \cite[Theorem 11.1.19(iv)]{Max19}. 
The polarization is independent of the  resolution by Lemma \ref{Lem intersection form}. 

\hfill
\qedsymbol

\end{Pf}

\subsection{Period maps}

Following the notation of \cite{Amb05}, we fix the following setup: 

\begin{itemize}
    \item $(X,B)$ is a klt pair of dimension $d$, such that $K_X + B \sim_{\mathbb{Q}} 0$. 

    \item $f: X \to S$ is a projective contraction to a nonsingular algebraic variety $S$. 

    \item $\mu: Y \to X$ is a resolution of singularities, $K_Y +B_Y = µ^*(K_X+B)$ is the log pullback, and $E$ is the support of the fractional part of $B_Y$. We assume that $E$ has snc support and $µ_*T_Y\langle -E\rangle $ is a reflexive sheaf (such a resolution exists by \cite[Lemma 1.1]{Amb05}). 

    \item Let $m$ be the index of $K_Y+B_Y$ over the generic point of $S$, $\varphi$ be a function on $Y$ such that $(\varphi) = m(K_Y+B_Y)$ over the generic point of $S$ and $\tau:\tilde Y\to Y$ be the index $1$ cover with respect to $\varphi$. Let $\nu:V\to\tilde Y$ be a resolution of $\tilde{Y}$. 

    \item The families $(Y,B_Y) \to S$ and $V \to S$ are relatively log smooth over an open subset $U$ of $S$. Let $\kappa_s: T_{S,s} \to H^1(Y_s,T_{Y_s}\langle -E_s\rangle )$, $\kappa_{V_s}: T_{S,s} \to H^1(V_s,T_{V_s})$ be the induced Kodaira-Spencer class. 
    
    \item Let $\sigma^1_s: T_{S,s} \to Hom(H^0(Y_s, \mathcal{O}_{Y_s}(\lceil -B_Y \rceil), Ext^1_{Y_s}(\Omega_{Y_s}^1\langle E_s\rangle, \mathcal{O}_{Y_s}(\lceil -B_Y \rceil)$, $\sigma_{V_s}: Hom(H^{d,0}(V_s), H^{d-1,1}(V_s))$ be the maps induced by cup product with $\kappa_s$, $\kappa_{V_s}$, respectively. 

    \item Let $\Phi$ be the period map associated to the polarized variation of Hodge structure $(R^dh_*\mathbb{Z}_V)_{prim} \otimes_{\mathbb{Z}} \mathcal{O}_S$, where $h$ is the natural map $V\to S$. Let $\textnormal{d}\Phi$ be the differential of $\Phi$. 
\end{itemize}

We may summarize these data in the following commutative diagram.  

\[\begin{tikzcd}
	{(Y,B_Y)} & {\tilde{Y}} & V \\
	{(X,B)} \\
	S
	\arrow["\mu"', from=1-1, to=2-1]
	\arrow["\tau"', from=1-2, to=1-1]
	\arrow["\nu"', from=1-3, to=1-2]
	\arrow["h", from=1-3, to=3-1]
	\arrow["f"', from=2-1, to=3-1]
\end{tikzcd}\]

We need the following proposition from \cite{Amb05}: 

\begin{Prop}\label{Prop the same kernel}
\cite[Proposition 2.1]{Amb05} The maps $\kappa_s, \kappa_{V_s}, \sigma^1_s, \sigma_{V_s}, \textnormal{d}\Phi_s$ have the same kernel. 
\end{Prop}

From now on, we restrict to the following threefold setting. As explained after Theorem \ref{Thm family gives polarized variation of Hs}, once $H^1(V_s,\mathcal O_{V_s})=0$ is established, the whole $H^3(V_s)$ is primitive and its polarization, given by the intersection pairing, is independent of an ample divisor: 

\begin{itemize}
    \item Assume $f: (X,B)\to S$ is flat of relative dimension $3$, and a dense set of closed fibers are isomorphic to a fixed $\mathbb{Q}$-factorial strict terminal Calabi-Yau pair $(X_0,\Delta_0)$. So $(X,B)$ is terminal $\mathbb{Q}$-factorial by Lemma \ref{Lem deformation of divisor}. 
\end{itemize}

Let $\tilde{X}$ be the normalization of $X$ in $\mathbb{C}(\tilde{Y})$, and let $V^{min}$ be a minimal model of $V$ over $\tilde{X}$. Possibly replacing $V$ by a higher model, we may assume $V$ is a common resolution of $\tilde{Y}$ and $V^{min}$. 
We write $r:V\to V^{min}$ for the resolution and
$h^{min}:V^{min}\to\tilde X\to X\to S$ for the composite, retaining $h=h^{min}\circ r:V\to S$. 

\begin{Rem}
We briefly explain why we introduce $V^{min}$. 

The role of $V^{min}$ is to provide terminal models over the fixed normalized cyclic cover. Terminality  permits the use of Lemma \ref{Lem cohomology of Calabi-Yau type} and Lemma \ref{Lemma independence of resolution}. These inputs cannot be applied for $\tilde X_s$, which need not be terminal. Moreover, $X_s$ itself may have nontrivial canonical index, so its third cohomology need not contain the one-dimensional top Hodge piece supplied by the cyclic cover. 
\end{Rem}

We may organize them in the following diagram: 

\[\begin{tikzcd}
	& V \\
	{(Y,B_Y)} & {\tilde{Y}} & {V^{min}} \\
	{(X,B)} & {\tilde{X}} \\
	S
	\arrow[from=1-2, to=2-2]
	\arrow[from=1-2, to=2-3]
	\arrow[from=2-1, to=3-1]
	\arrow[from=2-2, to=2-1]
	\arrow[from=2-2, to=3-2]
	\arrow[from=2-3, to=3-2]
	\arrow[from=3-1, to=4-1]
	\arrow[from=3-2, to=3-1]
\end{tikzcd}\]

\begin{Lem}\label{Lem V_0}
A dense set of closed fibers of $V^{min} \to S$ is isomorphic over $\mathbb{C}$ to a fixed variety $V_0$ such that $V_0$ satisfies the conditions of Lemma \ref{Lem cohomology of Calabi-Yau type}. 
\end{Lem}

\begin{Pf}
Let $s\in U$ be a closed point such that $(X_s,B_s) \simeq_{\mathbb{C}}(X_0,\Delta_0)$. By \cite[Section 3]{HMX18a}, $V_s^{min}$ is a terminal variety such that $K_{V_s^{min}}$ is nef over $X_s$. 
We claim that $V_s^{min}$ has only finitely many possibilities, so one of them will be our $V_0$. 

Indeed, restricting to the fiber over $s$, we have the following diagram: 

\[\begin{tikzcd}
	& V_s \\
	{(Y_s,B_{Y_s})} & {\tilde{Y}_s} & {V^{min}_s} \\
	{(X_s,B_s)} & {\tilde{X}_s} \\
	\{s\}
	\arrow[from=1-2, to=2-2]
	\arrow[from=1-2, to=2-3]
	\arrow[from=2-1, to=3-1]
	\arrow[from=2-2, to=2-1]
	\arrow[from=2-2, to=3-2]
	\arrow[from=2-3, to=3-2]
	\arrow[from=3-1, to=4-1]
	\arrow[from=3-2, to=3-1]
\end{tikzcd}\]

By Lemma \ref{Lem restriction of index $1$ cover}, possibly shrinking $U$, ${\tilde{Y}_s}$ is isomorphic to the index $1$ cover of $(Y_s,B_{Y_s})$ over $Y_s$. 

Although $(Y_s,B_{Y_s})$ is non-unique, the function field extension $\mathbb{C}({\tilde{Y}_s})/\mathbb{C}(X_s)$ is uniquely determined by $(X_s,\Delta_s)$. By \cite[Proposition 2.1]{PW17}, possibly shrinking $U$ again, we may assume $\tilde{X}_s$ is normal. So $\tilde{X}_s$ is the normalization of $X_s$ in $\mathbb{C}({\tilde{Y}_s})$, which is uniquely determined by $(X_s,\Delta_s)$.

Let $V_s^{min,\prime}$ be another possibility obtained as the above process. Then the natural birational map $V_s^{min} \dashrightarrow V_s^{min,\prime}$ is crepant and small by \cite[Theorem 3.52, Corollary 3.54]{KM98}. 
Then in the sense of \cite[Definition 3.6.6]{BCHM}, $V_s^{min,\prime}$ is a weak log canonical model of $V_s^{min}$ over $\tilde{X}_s$, so $V_s^{min,\prime}$ has only finitely many possibilities by \cite[Theorem E]{BCHM}. 

Moreover, we show that all possible $V_s^{min}$ satisfy the conditions in Lemma $\ref{Lem cohomology of Calabi-Yau type}$, i.e. $V_s^{min}$ is terminal and $h^0(V_s^{min},\omega_{V_s^{min}}) = 1$. Since $V_s^{min}$ has terminal singularities, it suffices to show $h^0(V_s,\omega_{V_s}) = 1$. 

Since $(X_s,B_s)$ is a strict Calabi-Yau pair, we may run $K_{X_s}+(1+\varepsilon)B_s$-MMP and terminate with a klt Calabi-Yau variety $X_s^{min}$. Let $\tilde{X}_s^{min}$ be the index $1$ cover of $X_s^{min}$, then $\tilde{X}_s^{min}$ is birational to $V_s$. Pick a common resolution $W$ of $\tilde{X}_s^{min}$ and $V_s$, then $h^0(V_s,\omega_{V_s}) = h^0(W,\omega_{W}) = h^0(\tilde{X}_s^{min},\omega_{\tilde{X}_s^{min}}) = 1$. 

\hfill
\qedsymbol
\end{Pf}

\begin{Lem}\label{Lem restriction of index $1$ cover}
Let $s\in U$ be a closed point such that $(X_s,B_s) \simeq_{\mathbb{C}}(X_0,\Delta_0)$. Then possibly shrinking $U$, $\tilde{Y}_s$ is isomorphic to the index $1$ cover of $(Y_s,B_{Y_s})$ over $Y_s$. 
\end{Lem}

\begin{Pf}
Let $\hat{Y}_s$ be the index $1$ cover of $(Y_s,B_{Y_s})$. We note that the index $1$ cover is unique since we are working over $\mathbb{C}$, see,  for example,  \cite[Definition 5.19]{KM98}. 

By the assumption that $(Y,B_Y) \to S$ is log smooth over $U$, $(Y_s,B_{Y_s})$ is still a log resolution. In particular, for $s\in U$, $mB_Y$ is integral iff $mB_{Y_s}$ is integral. Recall that $m$ is the index of $K_Y+B_Y$ over the generic point $\eta$ of $S$, $\varphi$ is a function on $Y$ such that $(\varphi) = m(K_Y+B_Y)$ over $\eta$ and $\tilde{Y}$ is the index $1$ cover with respect to $\varphi$.

We may shrink $U$ so that $(\varphi_s)=m(K_{Y_s}+B_{Y_s})$ for every $s\in U$. For each $1\leq n<m$ for which $nB_Y$ is integral, put
\[
\mathcal L_n=\mathcal O_{Y_U}(n(K_{Y_U/U}+B_Y|_{Y_U})).
\]
Since $Y_U\to U$ is smooth and $nB_Y$ is integral, $\mathcal L_n$ is invertible, hence flat over $U$, and its restriction to each fibre is the line bundle  $\mathcal{O}_{Y_s}(n(K_{Y_s}+B_{Y_s}))$.
By minimality of $m$, its generic fibre has no nonzero sections. 
By upper semicontinuity, \cite[Tag 0BDN]{stacks-project}, each locus
\[
Z_n=\{s\in U:h^0(Y_s,\mathcal L_n|_{Y_s})\geq1\}
\]
is a proper closed subset. We may remove their finite union. So every remaining fibre has index $m$.

The cover $\tilde{Y}$ is given by the normalization of $\oplus_{i=0}^{m-1} \mathcal{O}_Y(\lfloor i(K_Y+B_Y)\rfloor)$, and the multiplication law is given by $\varphi$. Similarly, $\hat{Y}_s$ is given by the normalization of $\oplus_{i=0}^{m-1} \mathcal{O}_{Y_s}(\lfloor i(K_{Y_s}+B_{Y_s})\rfloor)$, and the multiplication law is given by $\varphi_s$. By \cite[Proposition 2.1]{PW17}, possibly shrinking $U$, we may assume $\tilde{Y}_s$ is normal. By uniqueness of normalization, we have $\tilde{Y}_s \simeq \hat{Y}_s$. 

So our result holds. 

\hfill
\qedsymbol
\end{Pf}

By stratification theory, see \cite[Corollary 8.1.15]{CEGT14}, we may shrink $U$ such that the higher direct image sheaf $R^3h^{min}_*\mathbb{Z}_{V^{min}}$ is locally constant on $U$. The sheaf $R^3h_*\mathbb{Z}_V$ is locally constant on $U$ since $V\to S$ is smooth over $U$, hence a locally trivial fibration over $U$ by \cite[Theorem 7.2.1]{CEGT14}. We note that locally constant property here is considered in \emph{the complex analytic topology}. 

To simplify notations, let $\mathcal{W}_{\mathbb{Z}}:= (R^3h_*\mathbb{Z}_V)/\{\text{Torsion}\}$, $\mathcal{W}: = \mathcal{W}_{\mathbb{Z}}\otimes_{\mathbb{Z}} \mathcal{O}_{U}$, $\mathcal{V}_{\mathbb{Z}}:= Im(R^3h^{min}_*\mathbb{Z}_{V^{min}} \to R^3h_*\mathbb{Z}_V)/\{\text{Torsion}\}$, $\mathcal{V}: = \mathcal{V}_{\mathbb{Z}}\otimes_{\mathbb{Z}} \mathcal{O}_{U}$. 

For a fibre with $(X_s,B_s)\simeq_{\mathbb C}(X_0,\Delta_0)$, let $Z_s$ be the minimal model appearing in the definition of strictness, and let $\widehat Z_s\to Z_s$ be its index-one cover. This cover is quasi-\'etale, so $\widetilde q(Z_s)=0$ implies $H^1(\widehat Z_s,\mathcal O_{\widehat Z_s})=0$. The variety $\widehat Z_s$ has rational singularities and is birational to $V_s$. A common resolution therefore gives $H^1(V_s,\mathcal O_{V_s})=0$.

After restricting to a connected smooth locus of $h$, the Hodge numbers of its fibres are locally constant. Hence $H^1(V_s,\mathcal O_{V_s})=0$ for every $s\in U$. 
By proper base change theorem, we have $R^1h_*\mathbb{Z}_V = 0$. 
Theorem \ref{Thm family gives polarized variation of Hs} therefore equips $\mathcal W_{\mathbb Z}$ with the polarized variation $(\mathcal W_{\mathbb Z},\mathcal F^\bullet_{\mathcal W},Q_{\mathcal W})$, with $Q_{\mathcal W}$ the intersection pairing. Note that $(\mathcal W_{\mathbb Z},\mathcal F^\bullet_{\mathcal W},Q_{\mathcal W})$ is independent of the choice of an ample divisor, as discussed in Remark \ref{Rem independence}. 

Since $V_s^{min}$ is terminal and $V_s\to V_s^{min}$ is a resolution, rationality of its singularities also gives
$H^1(V_s^{min},\mathcal O_{V_s^{min}})=0$. 

We claim that the sheaf $\mathcal{V}_{\mathbb{Z}}$ is naturally equipped with a polarized variation of Hodge structure $(\mathcal{V}_{\mathbb{Z}}, \mathcal{F}^\bullet_\mathcal{V}, Q_\mathcal{V})$: 

\begin{Lem}\label{Lem intersect pvhs}

We define  $\mathcal{F}^{p}_{\mathcal{V}} =  \mathcal{F}^{p}_{\mathcal{W}} \cap \mathcal{V} \subset \mathcal{V}$, and $Q_{\mathcal{V}}$ is defined by the restriction of $Q_{\mathcal{W}}$. 

Then $(\mathcal{V}_{\mathbb{Z}}, \mathcal{F}^\bullet_\mathcal{V}, Q_\mathcal{V})$ is a polarized variation of Hodge structure. 

\end{Lem}

\begin{Pf}
In the sense of Definition \ref{Def 7.3.4}, to show that $(\mathcal{V}_{\mathbb{Z}}, \mathcal{F}^\bullet_\mathcal{V})$ is a variation of Hodge structure, we need to verify: 

(0) The Hodge filtration $\mathcal{F}^{p}_{\mathcal{V}}$ is given by holomorphic subbundles; 

(1) $\mathcal{V} = \mathcal{F}^{p}_{\mathcal{V}} \oplus \overline{\mathcal{F}^{4-p}_{\mathcal{V}}}$ as $C^{\infty}$-vector bundles; 

(2) Griffiths transversality.

(0) It suffices to prove that the Hodge numbers of $Im(H^3(V^{min}_s) \to H^3(V_s))$ are constant. Indeed, since the natural map $H^3(V^{min}_s) \to H^3(V_s)$ is a morphism of mixed Hodge structures, $Im(H^3(V^{min}_s) \to H^3(V_s))$ carries a pure Hodge structure. Moreover, by Lemma \ref{Lem cohomology of Calabi-Yau type}, the $(3,0)$ part of $Im(H^3(V^{min}_s) \to H^3(V_s))$ is of dimension $1$, so is $(0,3)$-part. Then $(1,2)$-part  and $(2,1)$-part are of constant dimension for topological reasons. 

(1) Since $\mathcal{W} = \mathcal{F}^{p}_{\mathcal{W}} \oplus \overline{\mathcal{F}^{4-p}_{\mathcal{W}}}$ as $C^{\infty}$-vector bundles, we have $\mathcal{F}^{p}_{\mathcal{V}} \cap \overline{\mathcal{F}^{4-p}_{\mathcal{V}}} = 0$. But $\dim \mathcal{F}^{p}_{\mathcal{V}} + \dim \overline{\mathcal{F}^{4-p}_{\mathcal{V}}} = \dim \mathcal{V}$ holds pointwise, so we have $\mathcal{V} = \mathcal{F}^{p}_{\mathcal{V}} \oplus \overline{\mathcal{F}^{4-p}_{\mathcal{V}}}$ as $C^{\infty}$-vector bundles. 

(2) It follows from Griffiths transversality for $\mathcal{W}$. 

To show that $\mathcal{V}$ is compatible with $Q_{\mathcal{V}}$, we need to prove the Hodge-Riemann relation, which can be verified pointwise, so the result follows from Lemma \ref{Lemma independence of resolution}. 

\hfill
\qedsymbol
\end{Pf}

\begin{Rem}
Pointwise, $\mathcal{V}$ is given by Lemma \ref{Lemma independence of resolution}. 
\end{Rem}

Recall that $\Phi: U \to \Gamma \backslash \mathcal{D}$ is the period map associated to the polarized variation of Hodge structure $(\mathcal{W}_{\mathbb{Z}}, \mathcal{F}^\bullet_\mathcal{W}, Q_\mathcal{W})$. Let $\Phi^{min}: U \to \Gamma^{min} \backslash \mathcal{D}^{min}$ be the period map associated to the polarized variation of Hodge structure $(\mathcal{V}_{\mathbb{Z}}, \mathcal{F}^\bullet_\mathcal{V}, Q_\mathcal{V})$. Let $\Gamma_a^{min}$ be the arithmetic monodromy group of $\mathcal{D}^{min}$.

\begin{Prop}\label{Prop inclusion of kernel}
The kernel of $\textnormal{d}\Phi^{min}_s$ is contained in the kernel of $\textnormal{d}\Phi_s$. 
\end{Prop}

\begin{Pf}
To prove this, we may work on a disk $C$ around $s$, and we may ignore the monodromy action on the period domain. 

Let $\tilde{\Phi}: C \to \mathcal{D}$ be the local lifting of $\Phi$, and $\tilde{\Phi}^{min}: C \to \mathcal{D}^{min}$ be the local lifting of $\Phi^{min}$, whose existence is given by Theorem \ref{Thm existence of period map}. 
Let $v=\dim \mathcal{V}$, $w = \dim \mathcal{W}$.

By Lemma \ref{Lem cohomology of Calabi-Yau type}, we have
\[
\mathcal F^3_{\mathcal V}=\mathcal F^3_{\mathcal W},
\qquad \operatorname{rank}\mathcal F^3_{\mathcal V}=1.
\]

Then by forgetting the Hodge filtrations except $\mathcal{F}^3$, we have holomorphic maps $f: \mathcal{D} \to \mathbb{C}\mathbb{P}^{w-1}$ and $f^{min}: \mathcal{D}^{min} \to \mathbb{C}\mathbb{P}^{v-1}$. Furthermore, we have a commutative diagram: 
\[\begin{tikzcd}
	C & {\mathcal{D}} & {\mathbb{C}\mathbb{P}^{w-1}} \\
	& {\mathcal{D}^{min}} & {\mathbb{C}\mathbb{P}^{v-1}}
	\arrow["\Phi", from=1-1, to=1-2]
	\arrow["{\Phi^{min}}"', from=1-1, to=2-2]
	\arrow["f", from=1-2, to=1-3]
	\arrow["{f^{min}}"', from=2-2, to=2-3]
	\arrow["\iota"', from=2-3, to=1-3]
\end{tikzcd}\]
where $\iota$ is induced by the natural inclusion from $\mathcal{V}$ to $\mathcal{W}$. 

Let $v\in T_{S,s}$ be a tangent vector satisfying $\mathrm d\widetilde\Phi^{min}_s(v)=0$.
Then its top Hodge piece is stationary to first order, i.e. $\mathrm d(f^{min}\circ\widetilde\Phi^{min})_s(v)=0$. The constant inclusion $\iota$ makes the top Hodge piece stationary in $\mathcal W$, i.e. $\mathrm d(\iota\circ f^{min}\circ\widetilde\Phi^{min})_s(v)=0$, which by definition gives  $\sigma_{V_s}(v)=0$. 
By \ref{Prop the same kernel}, we have
\[
\ker(\mathrm d\widetilde\Phi_s)=\ker(\sigma_{V_s}).
\]
Consequently $\mathrm d\widetilde\Phi_s(v)=0$, so we conclude that 
$\ker(\mathrm d\Phi_s^{min})\subseteq\ker(\mathrm d\Phi_s)$.

We note that no equality of the $(2,1)$-pieces is required in this proof. Our  conclusion uses Proposition \ref{Prop the same kernel} as an essential input and this is not a general property of subvariations of Hodge structures.

\hfill
\qedsymbol

\end{Pf}

\begin{Lem}\label{Lemma must be a point}

In the special threefold setting above, suppose that the Kodaira--Spencer map
\[
\kappa_s:T_{S,s}\longrightarrow H^1(Y_s,T_{Y_s}\langle-E_s\rangle)
\]
from Proposition \ref{Prop the same kernel} is injective for general $s\in U$. Then $S$ is a point. 
\end{Lem}

\begin{Pf}
Possibly shrinking $U$, we may assume $U$ is affine. 

By Proposition \ref{Prop inclusion of kernel}, $\textnormal{d}\Phi^{min}_s$ is also injective for general $s\in U$. Let $p: \Gamma^{min}\backslash \mathcal{D}^{min} \to \Gamma_a^{min}\backslash \mathcal{D}^{min}$ be the natural projection. Consider the map $\Phi^{min}_a = p \circ \Phi^{min}: U \to \Gamma^{min}_a\backslash \mathcal{D}^{min}$, which extends to a proper analytic map $\overline{\Phi^{min}_a}: \bar{U} \to \Gamma^{min}_a\backslash \mathcal{D}^{min}$ by \cite[Theorem 9.5, Proposition 9.11]{Gri70III}, where $\bar{U}$ is quasi-projective and contains $U$ as an open subset. Note that \cite[Theorem 9.5, Proposition 9.11]{Gri70III} are only stated for the monodromy group $\Gamma$, but the case for $\Gamma_a$ is also true by \cite[Theorem 3.1]{Sze99}. 

Let $t_0$ be the period point in $\Gamma_a^{min}\backslash \mathcal{D}^{min}$ corresponding to the Hodge structure of $Im(H^3(V_0,\mathbb{Z}) \to H^3(W,\mathbb{Z}))$, where $W$ is a resolution of $V_0$. This $t_0$ is independent of the resolution $W$ by Lemma \ref{Lemma independence of resolution}. Here $\Gamma_a^{min}$ cannot be replaced by $\Gamma^{min}$ since $\Gamma_a^{min}$ parametrizes isomorphism classes of polarized Hodge structures. 

Then $\Phi^{min,-1}_a(t_0)$ is a compact analytic subset of $\bar{U}$, which is algebraic by Chow's lemma. But $\Phi^{min,-1}_a(t_0)$ is also Zariski dense by Lemma \ref{Lem V_0}, so $\Phi^{min,-1}_a(t_0) = \bar{U}$. 

This means that $\Phi^{min}_a$ is constant. On a connected sufficiently small analytic neighbourhood, a local lift to $\mathcal D^{min}$ has image contained in one $\Gamma_a^{min}$-orbit. The integral isometry group acts properly discontinuously, and this orbit is discrete. The continuous local lift is therefore constant. But the fact that $\textnormal{d}\Phi^{min}_s$ is injective shows that $U$ must be a point, so is $S$. 

\hfill
\qedsymbol
\end{Pf}

\begin{Prop}\label{Proposition splitting as a product}
In the special threefold setting of this subsection, let $f:(X,B)\to S$ be a flat projective family over a smooth irreducible base, with $(X,B)$ klt, $K_X+B\sim_{\mathbb Q,S}0$ and with a Zariski-dense set of fibres isomorphic over $\mathbb C$ to the fixed strict terminal $\mathbb Q$-factorial Calabi--Yau pair $(X_0,\Delta_0)$. Assume the resolutions, cyclic covers and relative minimal models have been chosen as above. There exists a dominant generically finite morphism $\tau_U:U\to S$ such that
\[
(X,B)\times_S U\simeq_{\mathbb C}(X_0,\Delta_0)\times_{\mathbb C}U
\]
over $U$. 
\end{Prop}

\begin{Pf}
By \cite[Theorem 2.2]{Amb05}, there exist dominant morphisms $\tau : \bar{S} \to S$ and $\rho : \bar{S} \to S^{!}$, with $\tau$ generically finite and $\bar{S}, S^{!}$ nonsingular, and there exists a log variety $(X^!,B^!)$ with $K_{X^{!}} + B^{!} \sim_{\mathbb{Q}} 0$ and a projective contraction $f^!: X^! \to S^!$
\[\begin{tikzcd}
	{(X,B)} && {(X^!,B^!)} \\
	S & {\bar{S}} & {S^!}
	\arrow[from=1-1, to=2-1]
	\arrow[from=1-3, to=2-3]
	\arrow[from=2-2, to=2-1]
	\arrow[from=2-2, to=2-3]
\end{tikzcd}\]
such that 

(1)  there exists an open dense subset $U \subset \bar{S}$ and an isomorphism
\[\begin{tikzcd}
	{(X,B)\times_{S}\bar{S}|_U} && {(X^!,B^!)\times_{S^!}\bar{S}|_U} \\
	& U
	\arrow["\simeq", from=1-1, to=1-3]
	\arrow[from=1-1, to=2-2]
	\arrow[from=1-3, to=2-2]
\end{tikzcd}\]

(2) $\kappa_{s^!}$ is injective for general points $s^! \in S^!$. 

By (2), Proposition \ref{Prop the same kernel}, and Lemma \ref{Lemma must be a point}, $S^!$ must be a point. So $(X,B)\times_{S}\bar{S}|_U$ is a product family. So we may take $\tau_U$ to be the composition of the inclusion $U \to \bar{S}$ and $\tau: \bar{S} \to S$. 

\hfill
\qedsymbol
    
\end{Pf}

\subsection{Loci of fibers isomorphic to a fixed Calabi--Yau pair}

The following is our main theorem in this section. It parametrizes fibres isomorphic to a fixed strict terminal $\mathbb Q$-factorial Calabi--Yau pair. 

\begin{Thm}\label{Thm extraction of strict Calabi-Yau pairs}
Let $\pi: (\mathcal{X},\mathcal{D}) \to \mathcal{S}$ be a projective contraction of  varieties, and $(X_0,\Delta_0)$ be a strict terminal $\mathbb{Q}$-factorial Calabi-Yau pair of dimension $3$. 

Then there exists a finite type scheme $\mathcal{T}$ and a finite type morphism $f: \mathcal{T} \to \mathcal{S}$ such that $(\mathcal{X}_\mathcal{T},\mathcal{D}_\mathcal{T})$ splits as a product, i.e. there exists an isomorphism $(\mathcal{X}_\mathcal{T},\mathcal{D}_\mathcal{T}) \simeq \mathcal{T} \times_{\mathbb{C}}(X_0,\Delta_0)$, and the image of $f$ contains all closed points $s\in \mathcal{S}$ such that $(\mathcal{X}_s,\mathcal{D}_s) \simeq (X_0,\Delta_0)$. 
\end{Thm}

\begin{Pf}
By passing to a flattening stratification, we may assume $\pi$ is flat. As in \cite[Proposition 2.10]{BDCS20}, by passing to a stratification and possibly discarding some components, we may assume that:

1) $\mathcal{S}$ is smooth,

2) $\mathcal{D}$ does not contain any fibre of $\pi$, 

3) every fibre $\mathcal{X}_s$ is a normal variety, and

4) for any $s \in \mathcal{S}$ and any irreducible component $\mathcal{D}^\prime$ of $\mathcal{D}$, $\mathcal{D}^\prime_s$ is either empty or an irreducible prime divisor on $\mathcal{X}_s$. 

Furthermore, by discarding some open subsets, we may assume the set of points $s\in \mathcal{S}$ such that $(\mathcal{X}_s,\mathcal{D}_s) \simeq (X_0,\Delta_0)$ is dense in $\mathcal{S}$. 

Since $(\mathcal{X}_s,\mathcal{D}_s)$ is terminal $\mathbb{Q}$-factorial for a dense subset of $\mathcal{S}$, by Lemma \ref{Lem deformation of divisor}, $(\mathcal{X},\mathcal{D})$ is terminal $\mathbb{Q}$-factorial. 
For any ample divisor $\mathcal{H}$ on $\mathcal{X}$ and any $n\in \mathbb{Z}$, $\mathcal{H}+n(K_{\mathcal{X}}+\mathcal{D})$ is ample over $\mathcal{S}$ since this holds for a dense set of fibers and ampleness is an open condition. 
In particular, $\mathcal{H}+n(K_{\mathcal{X}}+\mathcal{D})$ is ample over the generic point $ \eta_\mathcal{S}$ of $\mathcal{S}$. 
This means $K_{\mathcal{X}}+\mathcal{D}$ is numerically trivial over $\eta_\mathcal{S}$. By \cite[Corollary 1.6]{HX11}, $K_{\mathcal{X}}+\mathcal{D}$ is $\mathbb{Q}$-linearly trivial over $\eta_\mathcal{S}$. In particular, after shrinking $\mathcal{S}$, we may assume $K_{\mathcal{X}}+\mathcal{D}$ is $\mathbb{Q}$-linearly trivial over $\mathcal{S}$. 

Now $\pi$ is a klt-trivial fibration, and a dense set of closed fibers are isomorphic to $(X_0,\Delta_0)$. By Proposition \ref{Proposition splitting as a product}, there exists a generically finite base change $f: \mathcal{T} \to \mathcal{S}$ such that $\pi$ splits as a product. Then we may discard the image of $f$ and apply Noetherian induction. 

\hfill
\qedsymbol
\end{Pf}

\section{Finiteness of fiber space structures}

In this section, we prove finiteness of fiber space structures of strict klt Calabi-Yau pairs in dimension 3 by combining boundedness and discreteness. Our main result is Theorem \ref{Thm finiteness for klt}. 

We start with an interesting observation: two different surface fibrations are dominated by a curve fibration. 

\begin{Lem}\label{Lemma factorization}
Let $f_i:X\to C_i$, $i=1,2$ be two different surface fibrations of a normal threefold $X$. Then they are dominated by a curve fibration. 
\end{Lem}

\begin{Pf}
Consider $\pi = (f_1,f_2):X \to C_1 \times C_2$. 

We claim that $\pi$ is surjective. To prove this, possibly replacing $X$ with a log resolution, we may assume $X$ is smooth. Let $A_i$ be an ample divisor on $C_i$, $i=1,2$, $H$ be a general member of a sufficiently ample class on $X$. By Bertini theorem, $H$ is smooth. If $\pi$ is not surjective, we have $(f_1^*A_1+f_2^*A_2)^2.H = 0$, hence $f_1^*A_1.f_2^*A_2.H = 0$, i.e. $(f_1^*A_1)|_H.(f_2^*A_2)|_H = 0$. By Hodge index theorem, $(f_1^*A_1)|_H$ and $(f_2^*A_2)|_H$ must be numerically proportional. 
Weak Lefschetz gives an injection
$H^2(X,\mathbb Q)\hookrightarrow H^2(H,\mathbb Q)$.
Using the identification of rational divisor classes modulo numerical equivalence with their rational first Chern classes on smooth projective varieties, it follows that
$N^1(X)\to N^1(H)$ is injective. 
So $f_1^*A_1$ and $f_2^*A_2$ are numerically proportional, which contradicts the assumption that $f_1$ and $f_2$ are different. 

Then the result follows by taking the Stein factorization. 

\hfill
\qedsymbol
\end{Pf}

From the boundedness result, we can prove finiteness of $v$-marked curve fibrations from a fixed strict terminal $\mathbb{Q}$-factorial Calabi-Yau pair of dimension $3$. 

\begin{Thm}\label{Theorem Finite elliptic}
Let $(X,\Delta)$ be a strict terminal $\mathbb{Q}$-factorial Calabi-Yau pair of dimension $3$. Fix a positive number $v$. Then the set of $v$-marked curve fibrations starting from $X$ is finite up to $Aut(X,\Delta)$. 
\end{Thm}

\begin{Pf}
By Theorem \ref{Theorem Boundedness}, the set of $v$-marked curve fibrations of strict terminal $\mathbb{Q}$-factorial Calabi-Yau pairs of dimension $3$ is bounded. 

By definition, we have a pair $(\mathcal{X},\mathcal{D})$ and  quasi-projective varieties $\mathcal{B}, \mathcal{S}$ with a $\mathbb{Q}$-Cartier $\mathbb{Q}$-divisor $\mathcal{A}$ on $\mathcal{B}$
and a commutative diagram 
\[\begin{tikzcd}
	({\mathcal{X}}, {\mathcal{D}}) && {\mathcal{B}} \\
	& {\mathcal{S}}
	\arrow["f", from=1-1, to=1-3]
	\arrow["g"', from=1-1, to=2-2]
	\arrow["h", from=1-3, to=2-2]
\end{tikzcd}\]
of projective morphisms such that for any $v$-marked curve fibration of a terminal $\mathbb{Q}$-factorial Calabi-Yau pair $(X,\Delta,B,\pi,A) \in V$ there exists a closed point $s\in \mathcal{S}$ such that $(X,\Delta,B,\pi,A) \simeq (\mathcal{X}_s, \mathcal{D}_s, \mathcal{B}_s, f_s, [\mathcal{A}_s])$, where $[-]$ denotes the numerical equivalence class. 

By Theorem \ref{Thm extraction of strict Calabi-Yau pairs}, there exists a finite type morphism $f: \mathcal{T} \to \mathcal{S}$ such that $(\mathcal{X}_\mathcal{T},\mathcal{D}_\mathcal{T})$ splits as a product, and the image of $f$ contains all closed points $s\in \mathcal{S}$ such that $(\mathcal{X}_s,\mathcal{D}_s) \simeq (X_0,\Delta_0)$. Let $f_{\mathcal{T}}^*\mathcal{A}_{\mathcal{T}}$ be the pullback of the polarization $\mathcal{A}_{\mathcal{T}}$ on $\mathcal{B}_{\mathcal{T}}$ via $f_{\mathcal{T}}$. 

This means that $(f_{\mathcal{T}}^*\mathcal{A}_{\mathcal{T}})_t$ falls into finitely many algebraic equivalence classes, hence $\pi^*A$ falls into finitely many numerical equivalence classes. Moreover, the same numerical equivalence class of $\pi^*A$ corresponds to isomorphic $v$-marked curve fibrations of $X$. So all possible $v$-marked curve fibrations starting from $X$ are finite up to $Aut(X,\Delta)$. 

\hfill
\qedsymbol
\end{Pf}

\begin{Rem}
This is stronger than finiteness of curve fibrations, since there exists a universal constant $v_0$ such that
any curve fibration of strict terminal Calabi-Yau pair of dimension $3$ admits a structure of $v_0$-marked curve
fibration.
\end{Rem}

From finiteness of $v$-marked curve fibrations from a fixed strict terminal $\mathbb{Q}$-factorial Calabi-Yau pair of dimension $3$, we can deduce an interesting result about lifting of automorphisms: 

\begin{Thm}\label{Theorem lifting}
Let $(X,\Delta)$ be a strict terminal $\mathbb{Q}$-factorial Calabi-Yau pair of dimension $3$, and $\pi: (X,\Delta)\to B$ be a curve fibration. Then there exists a finite index subgroup $G(B)$ of $Aut(B)$ that lifts to $Aut(X,\Delta)$ up to the subgroup $Aut^0(B)=\mathcal{A}ut^0(B)(\mathbb{C})$, i.e. for any $g\in G(B)$, there exists $s\in Aut^0(B), g_X\in Aut(X,\Delta)$, such that the following diagram commutes: 
\[\begin{tikzcd}
	(X,\Delta) & (X,\Delta) \\
	B & B
	\arrow["{g_X}", from=1-1, to=1-2]
	\arrow["\pi", from=1-1, to=2-1]
	\arrow["\pi", from=1-2, to=2-2]
	\arrow["g\circ s", from=2-1, to=2-2]
\end{tikzcd}\]
\end{Thm}

\begin{Pf}
Pick a very ample divisor $A$ on $B$. Then $\mathcal{F} : =(X,\Delta,B,\pi,A)$ is a $v$-marked curve fibration for some $v$. Denote by $\mathfrak{C}(X,\Delta,v)$ the set of $v$-marked curve fibrations of $(X,\Delta)$. Then $\mathfrak{C}(X,\Delta,v)$ is finite by Theorem \ref{Theorem Finite elliptic}. 

Let $L$ be the subgroup of ${Aut}(B)$ consisting of elements that lifts to ${Aut}(X,\Delta)$. We write $K = {Aut}(B,[A])$. 

Then consider the double coset $K\backslash{Aut}(B)/L$. Observe that each element $KgL \in K\backslash{Aut}(B)/L$ gives an element $(X,\Delta,B,\pi,g_*A)$, and for any $k\in K$, $l\in L$, we have $(X,\Delta,B,\pi,g_*A) \simeq (X,\Delta,B,\pi,(kgl)_*A)$ by the following diagram
\[\begin{tikzcd}
	{(X,\Delta)} & {(X,\Delta)} & {(X,\Delta)} \\
	{(B,g_*A)} & {(B,(kg)_*A)} & {(B,(kgl)_*A)}
	\arrow["{\tilde{k}}", from=1-1, to=1-2]
	\arrow["\pi", from=1-1, to=2-1]
	\arrow["id", from=1-2, to=1-3]
	\arrow["\pi", from=1-2, to=2-2]
	\arrow["\pi", from=1-3, to=2-3]
	\arrow["k", from=2-1, to=2-2]
	\arrow["id", from=2-2, to=2-3]
\end{tikzcd}\]
where $\tilde{k}$ is any lift of $k$ to ${Aut}(X,\Delta)$. 

Moreover, different elements of the double coset $K\backslash{Aut}(B)/L$ give different elements in $\mathfrak{C}(X,\Delta,v)$. 

So $K\backslash {Aut}(B)/L$ is finite since $\mathfrak{C}(X,\Delta,v)$ is finite. 

By \cite[Remark 2.6]{MZ18}, we have $[K:{Aut}^0(B)] <\infty$. So we deduce that the double coset ${Aut}^0(B)\backslash{Aut}(B)/L$ is finite. Since ${Aut}^0(B)$ is a normal subgroup of ${Aut}(B)$, we deduce that the coset ${Aut}(B)/({Aut}^0(B)\cdot L)$ is also finite. 

Now we may take $G(B) = {Aut}^0(B)\cdot L$. 

\hfill
\qedsymbol
\end{Pf}

We now prove finiteness of fiber space structures for strict klt pairs.

\begin{Thm}\label{Thm finiteness for klt}\label{Thm finiteness for terminal}
Let $(X,\Delta)$ be a strict klt Calabi-Yau pair of dimension $3$.
Then the set of fibrations starting from $X$ is finite up to $Aut(X,\Delta)$.
\end{Thm}

\begin{Pf}

\noindent\textbf{Step 1: the terminal $\mathbb{Q}$-factorial case}

Assume first that $(X,\Delta)$ is terminal and $\mathbb{Q}$-factorial.
By Theorem \ref{Theorem Finite elliptic} and the fact that there exists a universal constant $v_0$ such that
any curve fibration of a strict terminal $\mathbb{Q}$-factorial Calabi-Yau pair of dimension $3$ admits a structure
of $v_0$-marked curve fibration, the set of curve fibrations is finite up to $Aut(X,\Delta)$. 

So we only need to deal with surface fibration. If there is at most one surface fibration, the result follows immediately. 

Otherwise, there are at least two different surface fibrations starting from $X$. Let $\pi_1:(X,\Delta) \to B_1$, $\dots$, $\pi_n:(X,\Delta) \to B_n$ be all possible curve fibrations up to $Aut(X,\Delta)$. 

Note that by \cite[Theorem 5.4]{BDCS20}, each $B_i$ can be equipped with a boundary $\Delta_i$ such that $(B_i,\Delta_i)$ is klt Calabi-Yau. By the cone conjecture of klt Calabi-Yau surfaces, see \cite[Theorem 4.1]{Tot09}, there are only finitely many contractions starting from $B_i$ up to $Aut(B_i)$. Since $G(B_i)$ is of finite index in $Aut(B_i)$, where $G(B_i)$ is given by Theorem \ref{Theorem lifting}, there are only finitely many contractions starting from $B_i$ up to $G(B_i)$, say $p_{i1}: B_i\to C_{i1}$, $\dots$, $p_{ik_i}: B_i\to C_{ik_i}$. 

We claim that all possible surface fibrations are among $p_{ij}\circ\pi_i$ for some $1\leq i\leq n$, $1\leq j\leq k_i$ up to $Aut(X,\Delta)$. 

Let $f:(X,\Delta)\to C$ be a surface fibration starting from $(X,\Delta)$. By Lemma \ref{Lemma factorization}, every surface fibration is factored through by a curve fibration. So $f$ is factored as $(X,\Delta) \stackrel{\pi}{\longrightarrow} B \stackrel{p}{\longrightarrow} C$. By assumption, $\pi$ is isomorphic to some $\pi_i$ up to $Aut(X,\Delta)$, and for this chosen $i$ and this chosen isomorphism, $p$ is isomorphic to some $p_{ij}$ up to $G(B_i)\subset Aut(B_i)$. So there exists a commutative diagram where $u,v,g,h$ are isomorphisms and $g \in G(B_i)$: 
\[\begin{tikzcd}
	(X,\Delta) & (X,\Delta) & (X,\Delta) \\
	B & {B_{i}} & {B_i} \\
	C & {C_{ij}} & {C_{ij}}
	\arrow["u", from=1-1, to=1-2]
	\arrow["\pi", from=1-1, to=2-1]
	\arrow["{\pi_i}", from=1-2, to=2-2]
    \arrow["{\pi_i}", from=1-3, to=2-3]
	\arrow["v", from=2-1, to=2-2]
	\arrow["p", from=2-1, to=3-1]
	\arrow["g", from=2-2, to=2-3]
	\arrow["{p_{ij}}", from=2-2, to=3-2]
	\arrow["{p_{ij}}", from=2-3, to=3-3]
	\arrow["h"', curve={height=18pt}, from=3-1, to=3-3]
\end{tikzcd}\]

By Theorem \ref{Theorem lifting}, there exists some $s\in Aut^0(B_i)$ such that $g^{-1} \circ s^{-1}$ lifts to an automorphism of $(X,\Delta)$. Moreover, by Blanchard's Lemma, see \cite[Proposition 4.2.1]{BSU13}, there exists some $t\in Aut(C_{ij})$ such that $p_{ij} \circ s = t \circ p_{ij}$. In sum, we have the following commutative diagram, where the horizental maps are isomorphisms: 
\[\begin{tikzcd}
	(X,\Delta) & (X,\Delta) & (X,\Delta) & (X,\Delta) \\
	B & {B_{i}} & {B_i} & {B_i} \\
	C & {C_{ij}} & {C_{ij}} & {C_{ij}}
	\arrow["u", from=1-1, to=1-2]
	\arrow["\pi", from=1-1, to=2-1]
	\arrow["i", curve={height=-18pt}, from=1-2, to=1-4]
	\arrow["{\pi_i}", from=1-2, to=2-2]
	\arrow["{\pi_i}", from=1-3, to=2-3]
	\arrow["{\pi_i}", from=1-4, to=2-4]
	\arrow["v", from=2-1, to=2-2]
	\arrow["p", from=2-1, to=3-1]
	\arrow["g", from=2-2, to=2-3]
	\arrow["{p_{ij}}", from=2-2, to=3-2]
	\arrow["s", from=2-3, to=2-4]
	\arrow["{p_{ij}}", from=2-3, to=3-3]
	\arrow["{p_{ij}}", from=2-4, to=3-4]
	\arrow["h"', curve={height=18pt}, from=3-1, to=3-3]
	\arrow["t", from=3-3, to=3-4]
\end{tikzcd}\]
Hence $f:(X,\Delta)\to C$ appears as $p_{ij}\circ \pi_i$.

\medskip
\noindent\textbf{Step 2: the general strict klt case.}

Let $f:(Y,\Delta_Y)\to(X,\Delta)$ be a $\mathbb{Q}$-factorial terminalization.
The pair $(Y,\Delta_Y)$ is again a strict terminal $\mathbb{Q}$-factorial Calabi-Yau pair.

Consider the set of fiber space structures starting from $(X,\Delta)$, say $\{g_i: (X,\Delta) \to B_i\}_{i\in I}$. Then $g_i\circ f: (Y,\Delta_Y) \to B_i$ are fiber space structures starting from $(Y,\Delta_Y)$. By Step 1, the morphisms $g_i\circ f$ have finitely many orbits under $Aut(Y,\Delta_Y)$, say we have representatives $h_1,\dots,h_n$, where $h_i: (Y,\Delta_Y) \to C_i$. 

Each $h_i$ is a terminal Calabi-Yau fiber space of relative dimension at most $2$, so there are only finitely many intermediate contractions up to $Aut(Y,\Delta_Y)$ by \cite[Theorem 1.4]{Li23} and \cite[Corollary 2]{Xu24}. In particular, there are finitely many intermediate contractions of $h_i$ such that the factorization is given by finitely many birational contractions up to $Aut(Y,\Delta_Y)$, say $h_{ij}: (Y,\Delta_Y) \xrightarrow{f_{ij}} (X_{ij},\Delta_{ij}) \xrightarrow{g_{ij}} C_i$.

Choose those representatives for which $(X_{ij},\Delta_{ij})$ is isomorphic
to $(X,\Delta)$, and fix pair isomorphisms
$a_{ij}:(X,\Delta)\xrightarrow{\sim}(X_{ij},\Delta_{ij})$.
For any original $g:(X,\Delta)\to B$, the above reductions give an index $(i,j)$, a pair isomorphism $a: (X,\Delta) \to (X_{ij},\Delta_{ij})$, and a base isomorphisms $b: B\to C_{i}$
with $b\circ g=g_{ij}\circ a$.
Thus $a=a_{ij}\circ u$ for $u=a_{ij}^{-1}\circ a\in Aut(X,\Delta)$, and
\[
b\circ g=(g_{ij}\circ a_{ij})\circ u.
\]
Then the finitely many maps $g_{ij}\circ a_{ij}$ represent all the
required orbits. 

\hfill
\qedsymbol
\end{Pf}

\section{The cone conjecture for Calabi-Yau threefolds with positive augmented irregularity}

In this section, we prove the cone conjecture for klt Calabi-Yau pairs with positive augmented irregularity by using the Albanese map. Our main result is Theorem \ref{Thm cone of irregular}.

\subsection{The cone conjecture for quotients of abelian varieties}

We need the following result from \cite{MQ24}, with some slight modifications. 

\begin{Thm}\label{Thm cone of quotient of abelian varieties}
Let $A$ be an abelian variety, and $G$ be a finite group acting faithfully on $A$. Let $(X,\Delta)$ be the quotient of $A$ by $G$, where $\Delta$ is defined by ramification formula, see \cite[2.41]{Ko13}, then: 

(1) There exists a rational polyhedral cone $\Pi^\prime \subset (\bar{\mathcal{A}}(A)^G)_{+}$ such that $\mathcal{A}(A)^G \subset H\Pi^\prime$ for some subgroup $H$ of $N_{Aut(A)}(G)$, where $N_{Aut(A)}(G)$ is the normalizer of $G$ in $Aut(A)$. 

(2) The action of the image of the restriction map $N_{Aut(A)}(G) \to Aut(X,\Delta)$ on $\mathcal{A}(X)_+$ is of polyhedral type. 

(3) $\bar{\mathcal{A}}(X)_+ = \mathcal{A}^e(X) = \mathcal{B}^e(X)$. 

(4) The cone conjecture holds for $(X,\Delta)$. 
\end{Thm}

\begin{Pf}
(1) This follows from exactly the same arguments as \cite[Theorem A]{MQ24}. By \cite[Theorem 6.10]{MQ24}, $\mathcal{A}(A)^G$ is a homogeneous self-dual cone. By \cite[Proposition 6.15]{MQ24}, the image of $C_{Aut(A)}(G)$ is an arithmetic subgroup of $Aut(\mathcal{A}(A)^G)^0$. Applying \cite[Theorem 4.14]{MQ24}, we conclude (1). 

(2) Consider the following commutative diagram, where the right vertical map is an isomorphism by Lemma \ref{Lem cone of finite map}(1): 
\[\begin{tikzcd}
	{N_{Aut(A)}(G)} & {GL(N^1(A)^G)} \\
	{Aut(X,\Delta)} & {GL(N^1(X))}
	\arrow[from=1-1, to=1-2]
	\arrow[from=1-1, to=2-1]
	\arrow["\simeq", from=1-2, to=2-2]
	\arrow[from=2-1, to=2-2]
\end{tikzcd}\]

Combining Lemma \ref{Lem cone of finite map}(2), we conclude our result. 

(3) By \cite[Lemma 1.1]{Bau98}, a divisor $D$ on $A$ is algebraically equivalent to an effective divisor iff $D$ is nef. So $\bar{\mathcal{A}}(A)_+ = \mathcal{B}^e(A)$. By \cite[Lemma 5.1(1)]{LZ22}, $\mathcal{A}^e(A) \subset \bar{\mathcal{A}}(A)_+$. On the other hand, $ \bar{\mathcal{A}}(A)_+ = \bar{\mathcal{A}}(A)_+ \cap  \mathcal{B}^e(A) \subset \bar{\mathcal{A}}(A) \cap \mathcal{B}^e(A) = \mathcal{A}^e(A)$. So we have $\bar{\mathcal{A}}(A)_+ = \mathcal{A}^e(A) = \mathcal{B}^e(A)$. Then the result follows from Lemma \ref{Lem cone of finite map}(3)(5)(6).

(4) By (2) and Theorem \ref{Prop-Def}, there exists a rational polyhedral fundamental domain for the action of $Aut(X,\Delta)$ on $\mathcal{A}(X)_+$. So the cone conjecture holds by (3). 

\hfill
\qedsymbol

\end{Pf}

\subsection{Finish of the proof}

We need the following lemma for cones of surfaces: 

\begin{Lem}\label{Lem generalise in dim 2}
Let $f: (Y,0) \to (X,\Delta_X)$  be a finite Galois morphism of klt Calabi-Yau pairs of dimension at most $2$ with Galois group $G$, where the boundary divisors are defined by ramification formula. 

Let $H$ be the subgroup of $Aut(Y)$ that restricts to an automorphism of $(X,\Delta_X)$. 

Then the image of the descent map $r:H \to Aut(X,\Delta_X)$ is of finite index after taking the image to $GL(N^1(X))$. 
More precisely, let
$\rho_X:Aut(X,\Delta_X)\to GL(N^1(X))$ be the natural representation.
Then
\[
[\rho_X(Aut(X,\Delta_X)):\rho_X(r(H))]<\infty.
\]
\end{Lem}

\begin{Pf}
The $1$-dimensional case is trivial. We may only consider the surface case. 

\medskip
\noindent\textbf{Step 1: Reduce to the case when $Y$ is a canonical surface such that $K_Y\sim 0$. }

Consider the global index $1$ cover of $Y$, say $\pi:Z \to Y$. Then we claim that $f\circ \pi$ is Galois. 
By Lemma \ref{Lemma delicate lifting property} with zero boundary, any automorphism of $Y$ lifts to an automorphism of $Z$. Let $K$ be the subgroup of $Aut(Z)$ that restricts to an element in $G$. Then we have an exact sequence $0\to Aut(Z/Y) \to K \to G \to 0$, which implies $|Aut(Z/X)|\geq |K| \geq |Aut(Z/Y)|\cdot |G|$. So $f\circ \pi$ is Galois. 

We claim that if the result holds for $f\circ \pi$, then it also holds for $f$. Indeed, let $H^\prime$ be the subgroup of $Aut(Z)$ that restricts to an automorphism of $(X,\Delta_X)$. The group $H^\prime$ acts on intermediate fields of $\mathbb{C}(Z)/\mathbb{C}(X)$, which is a finite set. So a finite index subgroup $H^{\prime\prime}$ of $H^\prime$ fixes $\mathbb{C}(Y)$. By assumption, the image of the map $H^\prime \to Aut(X,\Delta_X)$ is of finite index after taking the image to $GL(N^1(X))$, so the image of the map $H^{\prime\prime} \to Aut(X,\Delta_X)$ is also of finite index after taking the image to $GL(N^1(X))$. But we have a factorization $H^{\prime\prime} \to H \to Aut(X,\Delta_X)$, so we conclude.

\medskip
\noindent\textbf{Step 2: Reduce from $Y$ to its smooth minimal resolution $W$.}

From the first step, we may assume $Y$ is a canonical surface such that $K_Y\sim 0$. 

Consider the minimal resolution of $Y$, say $\pi: W\to Y$. By uniqueness of minimal resolution, $G$ acts on $W$ and we may consider the quotient $Z = W/G$, with the boundary divisor $\Delta_Z$ defined by ramification formula. Then we have a commutative diagram: 
\[\begin{tikzcd}
	W & {Z=W/G} \\
	Y & {X = Y/G}
	\arrow["g", from=1-1, to=1-2]
	\arrow["\pi"', from=1-1, to=2-1]
	\arrow["\varpi", from=1-2, to=2-2]
	\arrow["f"', from=2-1, to=2-2]
\end{tikzcd}\]

In this case, $\pi$ and $\varpi$ are birational contractions. We claim that if the result holds for $g$, then it also holds for $f$. 

Suppose the result holds for $g$. Let $H^\prime$ be the subgroup of $Aut(W)$ that restricts to an automorphism of $(Z,\Delta_Z)$. Let $Aut(Z,\Delta_Z,\varpi)$ be the subgroup of $Aut(Z,\Delta_Z)$ that induces an automorphism of $(X,\Delta_X)$. 

By \cite[Lemma 3.4]{Tot09}, the action of $Aut(Z,\Delta_Z,\varpi)$ on $\mathcal{A}^e(X)$ admits a rational polyhedral fundamental domain. In view of \cite[Proposition 4.6]{Loo14}, the image of the map $Aut(Z,\Delta_Z,\varpi) \to Aut(X,\Delta_X)$ is of finite index after taking the image to $GL(N^1(X))$. 

By assumption, the image of the map $H^\prime \to Aut(Z,\Delta_Z)$ is of finite index after taking the image to $GL(N^1(Z))$. Let $H^{\prime\prime}$ be the subgroup of $Aut(W)$ that restricts to an element of $Aut(Z,\Delta_Z,\varpi)$. Then  the map $H^{\prime\prime} \to Aut(Z,\Delta_Z,\varpi)$ is of finite index after taking the image to $GL(N^1(Z))$.

Put $V=\varpi^*N^1(X)\subset N^1(Z)$. Then every element of $Aut(Z,\Delta_Z,\varpi)$ preserves $V$.
Restricting these linear maps to $V$, and identifying $V$ with $N^1(X)$
by $\varpi^*$, shows that the image of $H''$ has finite index in the image
of $Aut(Z,\Delta_Z,\varpi)$ on $N^1(X)$.

Combining the above results, we conclude that the image of the map $H^{\prime\prime} \to Aut(X,\Delta_X)$ is of finite index after taking the image to $GL(N^1(X))$. 

Note that elements in $H^{\prime\prime}$ induce automorphisms of $Y$, since the contracted curves are preserved. 
So we have a factorization $H^{\prime\prime} \to H \to Aut(X,\Delta_X)$, so the result holds.

\medskip
\noindent\textbf{Step 3:  Finish the proof.
}

After the reduction in Step 2, rename $W$ as $Y$ and $W/G$ as $X$.
Thus $Y$ in this step is a smooth K3 or abelian surface.

It suffices to show that $\mathcal{A}^e(X)$ admits a rational polyhedral fundamental domain under the action of $H$, by \cite[Proposition 4.6]{Loo14}. 

When $Y$ is a K3 surface, the result follows from \cite[Corollary 1.9]{OS01}. 

When $Y$ is an abelian surface, the result follows from Theorem \ref{Thm cone of quotient of abelian varieties}.

\hfill
\qedsymbol
\end{Pf}

\begin{Prop}\label{Prop cone conjecture of product type}
Let $S$ be a K3 surface with at worst canonical singularities, $E$ be an elliptic curve. Suppose we have a finite group $G$ acting on $S\times E$ diagonally, and the resulting quotient map $S\times E \to X: = (S\times E)/G$ is quasi-\'etale. 

Then the movable cone conjecture holds for $X$. 
\end{Prop}

\begin{Pf}

Denote $T = S/G$, $F=E/G$. Let $\Delta_T$ and $\Delta_F$ be the divisors defined by ramification formula. 

Consider the factorization $S\times E \stackrel{p}{\to} X= (S\times E)/G \stackrel{q}{\to} T\times F = (S/G)\times (E/G)$. By Lemma \ref{Lem cone of finite map}(1) and Corollary \ref{Cor cone of product}(1), we have $N^1(X) = N^1(S\times E)^G = (N^1(S)\times N^1(E))^G = N^1(S)^G\times N^1(E)^G = N^1(T)\times N^1(F) = N^1(T\times F)$. 
By Lemma \ref{Lem cone of finite map}(8) and Corollary \ref{Cor cone of product}(5), we have $\mathcal{M}^e(X) = \mathcal{M}^e(S\times E)^G = (\mathcal{M}^e(S)\times \mathcal{M}^e(E))^G = \mathcal{M}^e(S)^G\times \mathcal{M}^e(E)^G = \mathcal{M}^e(T)\times \mathcal{M}^e(F) = \mathcal{M}^e(T\times F)$. 

So $q^*$ induces isomorphisms $N^1(T\times F) \to N^1(X)$ and $\mathcal{M}^e(T\times F) \to \mathcal{M}^e(X)$.

By \cite[Theorem 4.1]{Tot09}, the cone conjecture holds for $(T,\Delta_T)$. 

Let $H$ be the subgroup of $Aut(S)$ that 
restricts to an automorphism of $(T,\Delta_T)$, i.e. the conjugation action of $H$ preserves the $G$-action on $S$. 
Let $\alpha: G \to Aut(S)$ be the action of $G$ on $S$. 
Then $H$ acts on $\alpha(G)$ via conjugation. Since $\alpha(G)$ is finite, the kernel of the action $H \to Aut(\alpha(G))$ is a finite index subgroup of $H$, say $H_0$. Let $H^\prime_0$ be the image of $H_0$ in $Aut(T,\Delta_T)$. Then by Lemma \ref{Lem generalise in dim 2}, the action of $H^\prime_0$ on $\mathcal{M}^e(T)$ is of polyhedral type. Let $\Pi^\prime$ be a rational polyhedral cone on $\mathcal{M}^e(T)$ such that $H^\prime_0\Pi^\prime = \mathcal{M}^e(T)$. Then $\Pi: = q^*(\Pi^\prime \times \mathcal{M}^e(F))$ is a rational polyhedral cone inside $\mathcal{M}^e(T\times F)=\mathcal{M}^e(X)$. 

Let $H_0$ act on $S\times E$ on the first factor. Then the conjugation action of $H_0$ commutes with the $G$-action on $S\times E$. So this induces an action of $H_0$ on $X$. Then we have $H_0\Pi = \mathcal{M}^e(X)$, so $Aut(X)\Pi = \mathcal{M}^e(X)$. By Theorem \ref{Prop-Def}, the cone conjecture holds for $X$.

\hfill
\qedsymbol
\end{Pf}

\begin{Rem}
The canonical divisor $K_X$ is numerically trivial by quasi-\'etale condition. 
\end{Rem}

Finally, we deal with the cone conjecture for  Calabi–Yau threefolds with positive
augmented irregularity. 

\begin{Thm}\label{Thm cone of irregular}
Let $X$ be a $\mathbb{Q}$-factorial klt Calabi-Yau threefold such that the augmented irregularity is positive. Then the movable cone conjecture holds for $X$. 
\end{Thm}

\begin{Pf}
Suppose that a finite quasi-\'etale cover $Z$ of $X$ has positive irregularity. Let $f: Z \to A$ be the Albanese map, by \cite[Theorem A]{Xu16},  there exists a finite \'etale cover $\pi: B\to A$ such that $Z\times_A B \simeq F\times B$. 

Let $\hat{F}$ be the index $1$ cover of $F$. By \cite[Theorem 1.5]{GKP16}, there exists a finite, quasi-\'etale Galois cover $\tilde{F}$ of $\hat{F}$ such that any finite \'etale cover of the regular locus $\tilde{F}^{0}$ of $\tilde{F}$ extends to a finite \'etale cover of $\tilde{F}$. 

When $\tilde{F}$ is an elliptic curve or an abelian surface, $X$ admits a quasi-\'etale cover, say $Z^\prime$, which is an abelian threefold. We may take the Galois closure, say $W \to Z^\prime \to X$. Then $W\to Z^\prime$ is quasi-\'etale, so it has to be \'etale by purity of ramification locus, see \cite[0EA1]{stacks-project}. But \'etale covers of abelian varieties are still abelian varieties, so $X$ is realized as a finite quasi-\'etale quotient of an abelian variety. Then the result follows from Theorem \ref{Thm cone of quotient of abelian varieties}. 

Otherwise, $\tilde{F}$ is a K3 surface with at worst canonical singularities. Then we have a finite quasi-\'etale map $h: \tilde{F}\times B \to X$. Then $\tilde{F}^{0}$ has no nontrivial finite \'etale cover since $\tilde{F}$ has no nontrivial finite \'etale cover. 

By taking the Galois closure, there exists a finite quasi-\'etale morphism $u: Y \to \tilde{F}\times B$ such that $Y/X$ is Galois. 
By purity of ramification locus, see \cite[0EA1]{stacks-project}, $u$ is finite \'etale when restricted to $u^{-1}(\tilde{F}^0\times B)$. By \cite[X. Corollaire 1.7]{SGA1}, $\pi_1^{\text{\'et}}(\tilde{F}^0\times B)$ is canonically isomorphic to $\pi_1^{\text{\'et}}(\tilde{F}^0)\times \pi_1^{\text{\'et}}(B) = \pi_1^{\text{\'et}}(B)$. This means the restriction of $u$ on $u^{-1}(\tilde{F}^0\times B)$ is given by a finite \'etale cover of $B$. 

In particular, there exists a finite \'etale morphism $v: C \to B$, such that $u|_{u^{-1}(\tilde{F}^0\times B)} = id_{\tilde{F}^0}\times v$. 
By uniqueness of normalization, $Y = \tilde{F} \times C$, $u = id_{\tilde{F}}\times v$. Then we have a finite Galois quasi-\'etale morphism $s = h \circ u: \tilde{F} \times C \to X$. Let $G$ be the Galois group. 

Then for $g\in G$, $g$ acts on $\tilde{F}\times C$ by $g(x,y) = (f_g(x,y),c_g(x,y))$. But there are no non-constant maps from $\tilde{F}$ to $C$, since $\tilde{F}$ has trivial irregularity and $C$ is an abelian variety. So $c_g(x,y) = c_g(y)$, and $g(x,y) = (f_g(x,y),c_g(y))$. This gives a map $C\to \mathcal{A}ut(\tilde{F})$, $y \mapsto f_g(-,y)$. But $\mathcal{A}ut(\tilde{F})$ is discrete when $\tilde{F}$ has trivial irregularity by \cite[Theorem 4.5]{Xu16}, so $f_g(-,y) = f_g(-)$ is independent of $y$. 

We conclude that $g(x,y) = (f_g(x),c_g(y))$. So $G$ acts on $\tilde{F}\times C$ diagonally. 

Then the result follows from Proposition \ref{Prop cone conjecture of product type}. 

\hfill
\qedsymbol
\end{Pf}

\section{Proof of main results}

\begin{Pf}(of Theorem \ref{Thm main})

(1) This is done by \cite[Corollary 2]{Xu24}. 

(2)(a) By Theorem \ref{Thm cone of irregular}, the movable cone conjecture holds for $X^{min}$. Then the movable cone conjecture holds for $(X,\Delta)$ by \cite[Theorem 2]{Xu24}. 

(2)(b) This is done by Theorem \ref{Thm finiteness for klt}. 

\hfill
\qedsymbol
\end{Pf}

\begin{Pf}(of Corollary \ref{Cor main})

When $X$ is $\mathbb{Q}$-factorial, this follows from \cite[2.5]{Xu24} and Theorem \ref{Thm main}. 

In general, let $f:(Y,\Delta_Y)\to(X,\Delta)$ be a terminal $\mathbb{Q}$-factorial modification of $(X,\Delta)$. 
Then we proceed as in Theorem \ref{Thm finiteness for klt}. 

Consider the set of fiber space structures starting from $(X,\Delta)$, say $\{g_i: (X,\Delta) \to B_i\}_{i\in I}$. Then $g_i\circ f: (Y,\Delta_Y) \to B_i$ are fiber space structures starting from $(Y,\Delta_Y)$. By assumption, $g_i\circ f$ is finite up to $Aut(Y,\Delta_Y)$, say we have representatives $h_1,\dots,h_n$, where $h_i: (Y,\Delta_Y) \to C_i$. 

Each $h_i$ are terminal Calabi-Yau fiber spaces of relative dimension at most $2$, so there are only finitely many intermediate contractions up to $Aut(Y,\Delta_Y)$ by \cite[Theorem 1.4]{Li23} and \cite[Corollary 2]{Xu24}. In particular, there are finitely many intermediate contractions of $h_i$ such that the factorization of given by birational contractions up to $Aut(Y,\Delta_Y)$, say $h_{ij}: (Y,\Delta_Y) \xrightarrow{f_{ij}} (X_{ij},\Delta_{ij}) \xrightarrow{g_{ij}} C_i$.

Choose those representatives for which $(X_{ij},\Delta_{ij})$ is isomorphic
to $(X,\Delta)$, and fix pair isomorphisms
$a_{ij}:(X,\Delta)\xrightarrow{\sim}(X_{ij},\Delta_{ij})$.
For any original $g:(X,\Delta)\to B$, the above reductions give an index $(i,j)$, a pair isomorphism $a: (X,\Delta) \to (X_{ij},\Delta_{ij})$, and a base isomorphisms $b: B\to C_{i}$
with $b\circ g=g_{ij}\circ a$.
Thus $a=a_{ij}\circ u$ for $u=a_{ij}^{-1}\circ a\in Aut(X,\Delta)$, and
\[
b\circ g=(g_{ij}\circ a_{ij})\circ u.
\]
Then the finitely many maps $g_{ij}\circ a_{ij}$ represent all the
required orbits. 

\hfill
\qedsymbol
\end{Pf}

\bibliographystyle{amsalpha}
\bibliography{refs}

\vspace{1em}

\noindent\small{\textsc{Qiuzhen College, Tsinghua University, Beijing, China} }

\noindent\small{Email: \texttt{xufl23@mails.tsinghua.edu.cn}}

\end{document}